\documentclass[11pt]{article} 
\usepackage{a4wide}
\usepackage{amsmath,amsfonts}

\newtheorem{theorem}{Theorem}

\newtheorem{proposition}[theorem]{Proposition}
\newtheorem{lemma}{Lemma} 
\newtheorem{claim}{Claim}

\numberwithin{lemma}{section}
\numberwithin{claim}{section}
\numberwithin{equation}{section}
 
\begin{document}

\title{Refined asymptotics around solitons for gKdV equations
 \footnote{
This research was supported in part by the Agence Nationale de la Recherche (ANR ONDENONLIN).
}}
\author{Yvan Martel$^{(1)}$ 
  and  Frank Merle$^{(2)}$ }
\date{\quad \\ (1)
 Universit\'e de Versailles Saint-Quentin-en-Yvelines,
Math\'ematiques \\
 45, av. des Etats-Unis,
 78035 Versailles cedex, France\\  
 martel@math.uvsq.fr\\
\quad  \\
(2)
Universit\'e de Cergy-Pontoise, IHES and CNRS, Math\'ematiques   \\
2, av. Adolphe Chauvin,
95302 Cergy-Pontoise cedex, France \\
 Frank.Merle@math.u-cergy.fr 
}
\maketitle
 
\begin{abstract}

We consider the generalized Korteweg-de Vries equation
\begin{equation}\label{kdvabs}
    \partial_t u + \partial_x (\partial_x^2 u + f(u))=0, \quad
    (t,x)\in [0,T)\times \mathbb{R}
\end{equation}
with general $C^2$ nonlinearity $f$. 
Under an explicit condition on $f$ and $c>0$,
there exists a solution  in the energy space $H^1$ of \eqref{kdvabs} of the  type $u(t,x)=Q_c(x-x_0-ct)$,  called soliton. Stability theory for $Q_c$ is well-known.

In  \cite{MMarchives},
\cite{MMnonlinearity}, we have proved that for $f(u)=u^p$, $p=2,3,4$,
the family of solitons 
is asymptotically stable in some local sense in $H^1$, i.e.
if $u(t)$ is close to $Q_{c}$ (for all $t\geq 0$), 
then $u(t,.+\rho(t))$ locally  converges in the energy space
to some $Q_{c_+}$  as $t\to +\infty$, for some $c^+\sim c$. The main improvement in \cite{MMnonlinearity}  is a direct proof, based on 
  a localized Viriel identity on the solution $u(t)$. As a consequence, we have obtained
an integral estimate on $u(t,.+\rho(t))-Q_{c_+}$ as $t\to +\infty$.

In \cite{yvanSIAM} and \cite{MMas1}, using the indirect approach of
\cite{MMarchives}, we could extend the asymptotic stability result
under general assumptions on $f$ and $Q_c$. However, without Viriel argument directly on the solution $u(t)$,
no integral estimate is available in that case.

The objective of this paper is twofold.

The main objective is to prove 
that in the case $f(u)=u^p$, $p=2,3,4$, $\rho(t)-c_+ t$ has limit as $t\to +\infty$
under the additional assumption $x_+ u\in L^2(\mathbb{R})$, which is consistent with a counterexample
in \cite{MMnonlinearity}. This result persists for general nonlinearity if a Virial type estimate is assumed. The main motivation for this type of result is the determination of explicit
shifts due to collision of two solitons in the nonintegrable case $p=4$, see \cite{MMcol1}.

The second objective of this paper is to provide large time stability and asymptotic stability
results for two soliton solutions for the case of general nonlinearity $f(u)$,
when the ratio of the speeds of the solitons is small.
The motivation is to accompany the two papers \cite{MMcol1}, \cite{MMcol2}, devoted to 
  collisions of two solitons in the nonintegrable case. 
The arguments are  refinements of \cite{We1}, \cite{MMT} specialized to the case
$u(t)\sim Q_{c_1}+Q_{c_2}$, for $0< c_2 \ll c_1$.
\end{abstract}

\section{Introduction}
We consider the generalized Korteweg-de Vries (gKdV) equations:
\begin{equation}\label{fkdv}
    \partial_t u + \partial_x (\partial_x^2 u + f(u))=0, \quad
    (t,x)\in [0,T)\times \mathbb{R},
\end{equation}
for $u(0)=u_0\in H^1(\mathbb{R})$,
with  a  general $C^2$ nonlinearity $f$.
We assume that 
\begin{equation}\label{surfas1}
\text{for $p=2$, $3$ or $4$,}\quad 
f(u)=  u^p+f_1(u)\quad \text{where  $\lim_{u\to 0}\left|\frac{f_1(u)}{u^p}\right|=0$.}
\end{equation}
  Denote $F(s)=\int_0^s f(s') ds'$. Note that for \eqref{fkdv}, one can solve the Cauchy locally in time
in $H^1$, using the arguments of Kenig, Ponce and Vega \cite{KPV}
(using the norms given for $f(u)=u^2$ in $H^1$). Moreover, the following conservation laws holds for $H^1$ solutions:
\begin{align*}
& \int u^2(t)=\int u_0^2,\quad E(u(t))=\frac 12 \int (\partial_x u(t))^2 -  \int F(u(t)) = \frac 12 \int (\partial_x u_{0})^2 -  \int F(u_0).
\end{align*}

Recall that if $Q_c$ is a solution of
\begin{equation}\label{ellipticas1} 
Q_c''+f(Q_c)=cQ_c, \quad x\in \mathbb{R}, \quad Q_c\in H^1(\mathbb{R}),
\end{equation}
then $R_{c,x_0}(t,x)=Q_c(x-x_0-ct)$
is solution of \eqref{fkdv}. We call soliton such nontrivial traveling wave solution
of \eqref{fkdv}. 

By well-known results on equation \eqref{ellipticas1} (see \cite{BL}, \cite{MMas1}), there exists
$c_*(f)>0$ such that
$$c_*(f)=\sup\{c>0 \text{ such that }
\forall c'\in (0,c), \text{ $\exists$ $Q_{c'}$ positive solution of 
\eqref{ellipticas1}}\}.$$
Note that  for $f(u)=u^p$, $c_*(u^p)=+\infty$ and for all $c>0$, 
$Q_c(x)=c^{\frac 1{p-1}} Q(\sqrt{c} x)$, where
$Q(x)=Q_1(x)=\left( \frac{p+1}2 {\cosh^{-2}\bigl(\frac {p-1} 2 \,x\bigr)}\right)^{\frac1 {p-1}}$.
Recall that if a solution $Q_c>0$ of \eqref{ellipticas1} 
exists then $Q_c$ is the unique (up to translation)
positive solution of \eqref{ellipticas1}
and  can be chosen  even on $\mathbb{R}$ and  decreasing on $\mathbb{R^+}$.

From Weinstein \cite{We1}, the soliton
$Q_c$ is orbitally stable if
\begin{equation}\label{stable}
\frac d{dc} \int {Q_{c'}^2(x) dx }_{\big |{c'}=c}>0.
\end{equation}

Concerning asymptotic stability, 
 we have proved in \cite{MMas1} the following general result

\smallskip

\noindent\textbf{Asymptotic stability (\cite{MMas1},\cite{MMnonlinearity},\cite{MMarchives})}\quad\textit{Assume that $f$ is $C^2$ and satisfies \eqref{surfas1}.
Let $0<c_0<c_*(f)$.
There exists $\alpha_0>0$ such that if
 $u(t)$ is a global  $(t\geq 0)$ $H^1$ solution of \eqref{fkdv}
satisfying}
\begin{equation}\label{th1bas1}
	\forall t\geq 0,\quad
	\inf_{r\in \mathbb{R}}\|u(t)-Q_{c_0}(.-r)\|_{H^1} < \alpha_0,
	\end{equation}
\textit{then the following hold.
\begin{enumerate}
\item Asymptotic stability in the energy space.
There exist  $t\mapsto c(t)\in (0,c_*(f))$, $t\mapsto \rho(t)\in \mathbb{R}$ such that
\begin{equation}
\label{th1-1}
u(t)-Q_{c(t)}(.-\rho(t))\to 0\quad \text{in $H^1(x>\tfrac {c_0 t}{10})$ as $t\to +\infty$.}
\end{equation}
\item Convergence of the scaling parameter.
Assume further that there exists $\sigma_0>0$ such that
$c\mapsto \int Q_c^2$ is not constant in any interval
$I\subset [c_0-\sigma_0,c_0+\sigma_0]$. 
Then, by possibly taking a smaller $\alpha_0>0$, 
there exits $c^+\in (0,c_*(f))$ such that $c(t)\to c^+$ as $t\to +\infty$.
Moreover, $\rho'(t)\to c^+$ as $t\to +\infty$.
\end{enumerate}}

\smallskip

Recall that the main improvement in \cite{MMnonlinearity} with respect to \cite{MMarchives} is a direct proof, based on  a localized Viriel etimate on the solution $u(t)$, see Claim \ref{VIR} in the present paper
for a similar result.
As a consequence, we have 
obtained the following   estimate: there exists $K>0$ such that
\begin{equation}
\label{th2A}
\int_0^{+\infty} \int (u(t,x)-Q_{c(t)}(x-\rho(t)))^2 e^{-\frac {\sqrt{c_0}} 4|x-\rho(t)|} dxdt \leq K \alpha_0^2.
\end{equation}
A Viriel identity was proved in the case $f(u)=u^p$ for $p=2,3,4$ 
see Proposition 6 in \cite{MMarchives},
and used under a localized form in \cite{MMnonlinearity} (see also the case $p=5$ in \cite{MMjmpa}).
In contrast, the proof of asymptotic stability in \cite{yvanSIAM},
\cite{MMas1} is for general $f(u)$, but it is indirect, following the original
approach of \cite{MMarchives}. Thus, in this case, it is not known whether \eqref{th2A} holds.
See further comments on this result in \cite{MMas1}.

\subsection{Refined asymptotics for   power nonlinearities}

Our main objective in this paper is to refine the convergence result for
the power case, i.e. $f(u)=u^p$ with $p=2,3$ and $4$ concerning the behavior of $\rho(t)$ as $t\to +\infty$.
In fact, the only requirement is that \eqref{th2A} holds, in particular, a virial type estimate around $Q_c$
is sufficient. A typical result in this direction is the following.

\begin{theorem}\label{TH2as2}
Assume $f(u)=u^p$, for $p=2,3$ or $4$.
Let $c_0>0$.
There exists $\alpha_0>0$ such that if
 $u(t)$ is an $H^1$ solution of \eqref{fkdv}
satisfying
\begin{equation}\label{th2AAA}
	\inf_{r\in \mathbb{R}}\|u(0,.+ r)-Q_{c_0}\|_{H^1} < \alpha_0,
\quad \int_{x>0} x^2u^2(0,x) dx <+\infty,
\end{equation}
then there exists
$c^+>0$ and $x^+\in \mathbb{R}$ such that
\begin{equation}
\label{th2B}
\lim_{t\to +\infty} c(t)=c^+,\quad 
\lim_{t\to +\infty} \rho(t)-c^+t=x^+,
\end{equation}
where $c(t)$ and $\rho(t)$ are defined as in \eqref{th1-1}.
\end{theorem}
Recall that Pego and Weinstein \cite{PW} and Mizumachi \cite{Mizu} 
also obtained   results of asymptotic
stability in weighted spaces, where convergence of $\rho(t)-t$ is proved.
The results \cite{PW} and \cite{Mizu} depend on a spectral assumption which is
proved only for $p=2,3$. Moreover, in \cite{PW},   the initial data  
has to belong to an exponential weigthed space. This condition has been relaxed in \cite{Mizu}, where the assumption is $\int_{x>0} x^{11} u^2 <+\infty$. 

\medskip

We point out two main motivations for this kind of results: 

$\bullet$ First, in \cite{MMnonlinearity}, we gave the following counterexample: 

\smallskip

\noindent \textit{For any $\alpha>0$,
there exists an $H^1$ solution $u(t)$ of the KdV equation i.e. (\ref{fkdv}) with $f(u)=u^2$, such that
$
\sup_{t\in \mathbb{R}} \|u(t)-Q(x-\rho(t))\|_{H^1}\le \alpha,
$ and for some $\kappa>0$, $\rho(t)$,}
\begin{equation*}
\lim_{t\rightarrow +\infty}\|u(t)-Q(x-\rho(t))\|_{H^1(x\ge t/2)}= 0,
\quad 
\lim_{t\rightarrow +\infty}\frac {\rho(t)-t}{\sqrt{\log(t)}} =\kappa.
\end{equation*}

\smallskip

The initial data used in this construction contains a series of small solitons, which converges in 
$H^1$ but not in $L^2(x^2dx)$. This proves that the convergence of $\rho(t)- t$ as $t\to +\infty$
is not true in general and 
requires some additional decay on the solution.  In this respect, the assumption $\int x^2u^2<+\infty$
in Theorem \ref{TH2as2} seems rather weak and improve the results in \cite{PW} and \cite{Mizu}.
When looking for an $L^2$ condition independent of $p$ for \eqref{fkdv}, we think that $\int x^2 u^2$
is optimal. Indeed, it seems that the relevant quantity is the $L^1$ norm, see \cite{MMcol1},
where it proved that during the collision of two solitons, the shifts on the trajectories are related
to $L^1$ norms.

Thus, the natural question left open by Theorem \ref{TH2as2} is whether assuming $\int_{x>0} |u|<+\infty$
is sufficient to obtain convergence of $\rho(t)-c^+ t$. This might require a much more refined analysis.

\medskip

$\bullet$
A second motivation for estimating $\rho(t)-c^+ t$ as $t\to \infty$ appears in the context of two soliton
collisions. In a paper \cite{MMcol1} 
concerning the collision of two solitons of different sizes for the gKdV equation
\eqref{fkdv} with $f(u)=u^4$, we were able to compute the shift on the trajectories of the solitons
resulting  from their collision. The explicit shift is obtained for a fixed time $T$, long after the
collision. The proof of Theorem \ref{TH2as2} allows us to quantify the variation of this shift due
to long time i.e. for  $t>T$.

\medskip

We also point out that the proof of Theorem \ref{TH2as2} relies on several refinements of 
the proof of the asymptotic stability in \cite{MMarchives}, \cite{MMnonlinearity} and
is independent from the methods in \cite{PW}, \cite{Mizu}. Let $\eta(t,x)=u(t,x)-Q_{c(t)}(x-\rho(t))$.
Observe
 that the standard estimate
$$|\rho'(t)-c(t)|\leq K \left(\int \eta^2(t,x) e^{-\frac{\sqrt{c_0}}4 |x-\rho(t)|} dx\right)^{\frac 12},$$
and \eqref{th2A} do not give any conclusion

We proceed as follows.
First, we improve the monotonicity arguments on $u(t)$ used in \cite{MMarchives}, \cite{MMT} and \cite{MMnonlinearity}. The improvement is to prove monotonicity results on $\eta(t)$, which are much more
precise (see Claim \ref{NEW}). This argument allows us to prove that \eqref{th2A} implies
$\int_0^{+\infty} |c(t)-c^+| dt <+\infty$. Second, the control of $\int_0^{+\infty} (\rho'(t)-c(t)) dt$
is obtained through the equation of $\eta(t)$, by noting that at the first order 
$\rho'(t)-c(t)$ is the derivative of some bounded function of $t$.
Note that we do not prove $\int_0^{+\infty} |\rho'(t)-c(t)| dt <+\infty$.

\subsection{Large time behavior in the two soliton case}
A second objective of this paper is to provide asymptotic analysis in large time related
to two soliton solutions of \eqref{fkdv}.

From  \cite{Martel} (see also \cite{MMT}), there exist solutions $u(t,x)$ of \eqref{fkdv}
which are asymptotic $N$-soliton solutions at $t\to -\infty$ in the following sense:
		\emph{let $N\geq 1$, $c_1>\ldots>c_N>0$, and $x_1,\ldots,x_N\in \mathbb{R}$, there exists a unique $H^1$ solution $U$  of (\ref{fkdv}) such that }
		\begin{equation}\label{limit0+}
            \lim_{t\to -\infty}
            \Bigl\|U(t)-\sum_{j=1}^N Q_{c_j }(.-x_j -c_j t)\Bigr\|_{H^1(\mathbb{R})}=0.
         \end{equation}

\medbreak

The behavior displayed by these solutions is stable in some sense.
Considering for example the case of two solitons,
there exist a large class of solutions such that as $t\sim -\infty$,
\begin{equation}\label{absun}
  u(t,x)=Q_{c_1}(x{-}x_1{-}c_1 t) + Q_{c_2}(x{-}x_2{-}c_2 t) + \eta(t,x),
\end{equation}
where $c_1>c_2$ and  $\eta(t)$ is a dispersion term small in the energy space $H^1$ with respect to $Q_{c_1}$, $Q_{c_2}$ (see \cite{MMT}).
From the Physics point of view, the two solitons $Q_{c_1}$ and $Q_{c_2}$ have collide at some time $t_0$.
In the special case
$c_2 \ll c_1$ (or equivalently, $\|Q_{c_2}\|_{H^1}\ll \|Q_{c_1}\|_{H^1}$) and  $\|\eta(t)\|_{H^1} \ll \|Q_{c_2}\|_{H^1}$, for $t$ close to $-\infty$,
we have introduced in \cite{MMcol1}
explicit computations allowing to understand the collision at the main orders,
using a new nonlinear ``basis'' 
to  write and compute an approximate solution $v(t,x)$ up to any order of size.

Recall that the problem of   collision of two solitons  is a classical 
question in nonlinear wave propagation. 
In the so-called integrable cases (i.e. $f(u)=u^2$ and $f(u)=u^3$) it is well-known that
there exist explicit multi-soliton solutions, describing the elastic collision of several solitons (see Hirota \cite{Hirota}, Lax \cite{Lax} and the review paper Miura \cite{Miura}).
Note that in experiments, or numerically for more accurate nonintegrable
models (see Craig et al. \cite{Craig}, Li and Sattinger \cite{LiSattinger}
and other references in \cite{MMcol1}),
this remarkable property is mainly preserved, i.e. the collision of two solitons is almost elastic, however, a (very small) residual part is observed after the collision.
Equation    \eqref{fkdv}  being not integrable (unless $f(u)=u^2$ and $f(u)=u^3$),
explicit $N$-soliton solutions are not  available in this case. 
The results obtained in \cite{MMcol1} and \cite{MMcol2}, using the present paper,
are the first rigorous results concerning inelastic (but almost elastic) collision in a nonintegrable situation.
We refer to the introduction  and
the references in \cite{MMcol1} for an overview on these questions.

In \cite{MMcol1} and \cite{MMcol2}, the approximate solution is adapted to treat a large but fixed time interval around the collision region, but not the large time asymptotics (see for example Proposition 3.1 in \cite{MMcol1}).
In Proposition \ref{TH3as2} below, we give stability and 
asymptotic stability results required to control the
asymptotics in large time of the solutions constructed in \cite{MMcol1} and \cite{MMcol2}. In particular, we need a sharp stability result in the case where
one soliton is small with respect of the other.  
We claim the following.
 \begin{proposition}\label{TH3as2}
Assume that $f$ is $C^2$ and satisfies \eqref{surfas1}
for $p=2$, $3$ or $4$. Let $0<c_2<c_1<c_*$ be such that \eqref{stable} holds for $c_1$, $c_2$.
Let $0<c_1<c_*(f)$ be such that \eqref{stable} holds. There exist  
$c_0(c_1)$ and $K_0(c_1)>0$, continuous in $c_1$
such that for any $0<c_2<c_0(c_1)$ and for any $\omega>0$, the following hold.
Let 
$T_{c_1,c_2}=c_1^{\frac 32} \big(\frac {c_2}{c_1}\big)^{-\frac 12-\frac 1{100}}.$
Let $u(t)$ be an $H^1$ solution of \eqref{fkdv} such that
for some $t_1\in \mathbb{R}$ and $\frac 12 {T_{c_1,c_2}}\leq X_0\leq \frac 32{T_{c_1,c_2}}$,
\begin{equation}
\label{D25intro}
\|u(t_1)-Q_{c_1}-Q_{c_2}(.+X_0)\|_{H^1}\le   c_2^{\omega+\frac 1{p-1}+\frac 14}.
\end{equation}
Then, there exist    $C^1$  functions $\rho_1(t)$, $\rho_2(t)$ defined on $[t_1,+\infty)$  such that
\begin{enumerate}
\item Stability of the two solitons.
\begin{equation}\label{huitintro}
\sup_{t\ge t_1} \| u(t)-Q_{c_1}(.-\rho_1(t))-Q_{c_2}(.-\rho_2(t))  \|_{H^1 }
\le K     c_2^{\omega+\frac 1{p-1}-\frac 14},
\end{equation}
\begin{equation}\label{supplintro}\begin{split}
&\forall t\ge t_1,\ \tfrac 12 c_1\leq (\rho_1-\rho_2)'(t) \leq \tfrac 32 c_1,
\\ &|\rho_1(t_1)|\leq K c_2^{\omega+\frac 1{p-1}+\frac 14},\quad
|\rho_2(t_1)-X_0|\leq K c_2^{\omega}.
\end{split}\end{equation}
\item Asymptotic stability.
There exist $c_1^+, c_2^+ >0$  such that
\begin{equation}\label{neufintro}
\lim_{t\rightarrow +\infty}\|u(t)-Q_{c_1^+}(x-\rho_1(t))-Q_{c_2^+}(x-\rho_2(t))\|
_{H^1(x> \frac{c_2t}{10})}=0.
\end{equation}
\begin{equation}\label{septintro}
\left|\frac{c_1^+}{c_1}-1\right|\leq   K    c_2^{\omega+\frac 1{p-1}+\frac 14 },\quad
 \left|\frac {c_2^+}{c_2} - 1\right|\le K  c_2^{\omega}.
\end{equation}
\end{enumerate}
\end{proposition}

\noindent\emph{Remark.} The time $T_{c_1,c_2}$ corresponds to a time long after the collision of the two 
solitons (see \cite{MMcol1}). In \eqref{D25intro}, since $X_0>\frac 12 T_{c_1,c_2}$,
the two solitons are decoupled for any $t\geq t_1$.

\medskip

Proposition \ref{TH3as2} follows directly from known arguments in Weinstein \cite{We1}, and
in \cite{MMT}, \cite{MMas1}. The only new point is the fact that one soliton is small with respect to the other.
However, these statements are essential in \cite{MMcol1} and \cite{MMcol2}. In those works, the point
is to show that even in a nonintegrable situation, the two soliton structure is preserved by collision.
The method in \cite{MMcol1}, \cite{MMcol2} concerns the collision problem 
in $[-T_{c_1,c_2},T_{c_1,c_2}]$. 
To obtain global in time results,
it is essential to prove the global in time stability of the two soliton structure after the collision, 
i.e. for $t>T_{c_1,c_2}$, which is provided by Proposition  \ref{TH3as2}.
Proposition \ref{TH3as2} is directly applied in \cite{MMcol2}.
A slightly more precise   stability result is required in \cite{MMcol1}, see 
Proposition \ref{ASYMPTOTIC} in Section 2.

\medskip

Now, we claim an extension of Theorem 1 to the case of two solitons, with a qualitative
control on $\lim_{t\to +\infty} \rho(t)-c^+ t$.

\begin{theorem}\label{TH4as2}
Under the assumptions of Proposition \ref{TH3as2}, 
assume further that $f(u)=u^p$, for $p=2,3$ or $4$ and
$\int_{x>0} x^2 \, u^2(0,x) dx <+\infty$. Then, there exist  $x_1^+$ and 
$x_2^+$ such that
\begin{equation}\label{decayfortintro}
\lim_{t\to +\infty} \rho_1(t)-c_1^+ t=x_1^+, \quad
\lim_{t\to +\infty} \rho_2(t)-c_2^+ t=x_2^+.
\end{equation}
In the case $p=4$, if in addition, for some $\kappa>0$,
\begin{equation}
\alpha< \kappa c_2^{\frac 13} \quad \text{and}\quad \int_{x>\frac {11} {12} |\ln c_2|} x^2 \, u^2(t_1,x) dx <  \kappa c_2^{\frac 54}
\end{equation}
then
\begin{equation}\label{decayfortintroquali}
|x_1^+-\rho_1(0)|\leq K c_2^{\frac 58}, \quad |x_2^+-\rho_2(0)|\leq K c_2^{\frac 1{12}}.
\end{equation}
\end{theorem}

The main motivation of Theorem \ref{TH4as2} is the following:
in \cite{MMcol1}, in the same context as before,
we were able to compute the main order of the shift on the trajectories of the solitons
due to the collision at time  $t=T_{c_1,c_2}$.
Theorem \ref{TH4as2} proves that the shifts  do   change at the main order in large time
(for example, at the main order, the shift of $Q_2$ is a nonzero constant independent
of $c_2$, so that it is preserved  by \eqref{decayfortintroquali}).
See proof of Theorem 1.2 in \cite{MMcol1} for details.

\medskip

The plan of the paper is as follows.
In Section 2, we prove the stability part of Proposition~\ref{TH3as2}. We focus on the case $f(u)=u^p$ for simplicity,  the proof in the general case being exactly the same (see
\cite{MMas1} and \cite{MMT}). Moreover, by a scaling argument we consider only the
case   $c_1=1$, $c_2=c$, where $c$ is small enough.

In Section 3, we prove   Theorem \ref{TH4as2}. First, we use the methods of
localized Viriel estimates as in \cite{MMnonlinearity} to obtain the equivalent of \eqref{th2A}
for two solitons. Next, we prove \eqref{decayfortintro} and \eqref{decayfortintroquali}.
The proof of Theorem \ref{TH2as2} follows directly from the arguments of Section 3, thus
it will be omitted.

\section{Stability for large time of  $2$-soliton like solutions}
Recall that we restrict ourselves to the case $f(u)=u^p$ ($p=2,3,4$) and $c_1=1$, $c_2=c$
small enough. Let 
$$T_c=c^{-\frac 12 -\frac 1{100}},\quad  q=\frac 1{p-1}-\frac 14.$$
Denote for $v\in H^1(\mathbb{R})$,
$
\|v\|_{H^1_c}=\left(\int_{\mathbb{R}} \left( (v'(x))^2 + c v^2(x) \right)dx \right)^{\frac 12},$
which corresponds to the natural norm to study the stability of $Q_c$.

\begin{proposition}[Stability of two decoupled solitons]\label{ASYMPTOTIC}
There exists $K>0$, $\alpha_0>0$, $c_0>0$  such that for any $0<c<c_0$, $0<\alpha<\alpha_0$,
the following is true.

 Let $u(t)$ be an $H^1$ solution of \eqref{fkdv} such that
for some $t_1\in \mathbb{R}$ and $X_0\geq  T_c/2$,
\begin{equation}
\label{D25}
\|u(t_1)-Q-Q_c(.+X_0)\|_{H^1}\le \alpha c^{q+\frac 12}.
\end{equation}
Then there exist    $C^1$  functions $\rho_1(t)$, $\rho_2(t)$ defined on $[t_1,+\infty)$  such that
\begin{equation}\label{huit}
\sup_{t\ge t_1} \| u(t)-(Q(.-\rho_1(t))+Q_c(.-\rho_2(t)))  \|_{H^1_c}
\le K \alpha  c^{q+\frac 12}+K\exp\left(-  {c^{-\frac 1{400}}}\right),
\end{equation}
\begin{equation}\label{suppl}\begin{split}
&\forall t\ge t_1,\ \tfrac 12\leq \rho_1'(t)-\rho_2'(t) \leq \tfrac 32,
\\ &|\rho_1(t_1)|\leq K \alpha c^{q+\frac 12},\quad
|\rho_2(t_1)-X_0|\leq K \alpha.
\end{split}\end{equation}
\end{proposition}
Proposition \ref{TH3as2} follows immediately from Proposition \ref{ASYMPTOTIC}
with $\alpha=c^\omega$ ($\omega>0$) and a scaling argument.

\medskip

\noindent\emph{Remark.}\quad Note that the proof of  Proposition \ref{ASYMPTOTIC} does not need any new arguments with respect to \cite{MMT}. We only need to check that the argument of \cite{MMT} still applies to the situation where one soliton is small with respect to the other.

Since $\|Q_c\|_{L^2}=c^q \|Q\|_{L^2}$ (see Claim \ref{LemmaA1}), the assumption \eqref{D25} does not seem optimal by a factor $\sqrt{c}$. This is due to the fact that the
appropriate norm for the stability of $Q_c$ is  $\|.\|_{H^1_c}$.

\medbreak

\noindent\emph{Proof of  Proposition \ref{ASYMPTOTIC}.} 
By time translation invariance, we may assume that $t_1=0$. Let  $X_0\ge T_c/2$ be such that
\begin{equation}\label{48bon}
\|u(0)-Q-Q_c(.+X_0)\|_{H^1}\le \alpha c^{q+\frac 12}.
\end{equation}

Let $D_0>2$  to be chosen later, $r=\frac 1{400}$  and 
\begin{eqnarray*}
t^*&=&\sup\Big\{t\geq 0 ~| ~\forall ~ t'\in [0,t), \ \exists \, \tilde\rho_1,\tilde\rho_2\in \mathbb{R} 
\ | \ \tilde \rho_1-\tilde \rho_2>\tfrac 14 T_c \\ &&\qquad
\text{and } \|u(t')-Q(.-\tilde\rho_1)-Q_c(.-\tilde\rho_2)\|_{H^1_c}\le D_0 (\alpha c^{q+\frac 12}+
\exp(-c^{-r}))\Big\}.
\end{eqnarray*}
Observe that $t^*>0$ is well-defined since $D_0>2$, (\ref{48bon}) and the continuity of $t\mapsto u(t)$ in ${H^1}$.
The objective is to prove $t^*=+\infty$. For the sake of contradiction, we assume that $t^*$ is finite.

First, we decompose the  solution on $[0,t^*]$ using modulation theory around the sum of two solitons (see proof of Claim \ref{LEMMEB1} in Appendix A.1).

\begin{claim}[Decomposition of the solution]\label{LEMMEB1}
For $\alpha>0$, $c>0$ small enough, independent of $t^*$,
there exist $C^1$  functions $\rho_1(t)$, $\rho_2(t)$, $c_1(t),$ $c_2(t)$,
defined on $[0,t^*]$,  such that the function $\eta(t)$ defined by
\begin{equation*}
\eta(t,x)=u(t,x)-R_1(t,x)-R_2(t,x),
\end{equation*}
where for $ j=1,2,$
$
R_j(t,x)=Q_{c_j(t)}(x-\rho_j(t)),$
satisfies for all $t\in [0,t^*],$
\begin{eqnarray}&&
\int R_j(t)\eta(t) =\int (x-\rho_j(t))  R_j(t)\eta(t)=0,\quad j=1,2,\label{dix}\\ &&
\|\eta(t)\|_{H^1}+|c_1(t)-1|+c^q\left|\frac {c_2(t)}c -1\right|\le K D_0 (\alpha  c^{q}+c^{-\frac 12}\exp(-c^{-r})),\label{onze}\\&&
|  \rho_2'(t)|+ |  \rho_1'(t)-1|\le {\textstyle\frac 1 {10}}, \quad  \rho_1(t)-\rho_2(t)\ge {\textstyle\frac t2}+ {\textstyle\frac 14}  T_c,\label{treize}\end{eqnarray}
\begin{equation}\label{quatorze}\begin{split}
&\|\eta(0)\|_{H^1}+|c_1(0)-1|+c^{q}
\left|\frac {c_2(0)}c - 1\right|\le K \alpha c^{q+\frac 12},
\quad  \rho_1(0)-\rho_2(0)\ge {\textstyle\frac 14} T_c,\\
&
|\rho_1(0)|+c^{q+\frac 12} |\rho_2(0)-X_0|\leq K \alpha c^{q+\frac 12}.
\end{split}\end{equation}
\end{claim}
We define 
\begin{equation}\label{surphiff} 
\begin{split}
&\psi(x)=\frac {2}{\pi}\arctan(\exp(- x/4)),
\text{ so that } {\rm lim}_{+\infty} \psi=0,\ {\rm lim}_{-\infty} \psi=1, \\ &
\forall x\in \mathbb{R},\quad  \psi(-x)=1-\psi(x),\quad  \psi'(x)=\frac 1{4 \pi {\rm cosh}(x/4)},
\quad  |\psi'''(x)|\le \frac 1{16} |\psi'(x)|.
\end{split}\end{equation}
For $m(t)=\frac 12 (\rho_1(t)+\rho_2(t))$, we set
\begin{equation}\label{defg}\begin{split}
& \mathcal{I}(t)=\int u^2(t) \psi\left(x-m(t)\right) dx,\quad 
 g(t)=\int \left(\eta_x^2(t,x) +  \left(c+ \psi(x-m(t))\right) \eta^2(t,x)\right)dx.
\end{split}\end{equation}
Note that $\mathcal{I}(t)$ corresponds at the main order to the $L^2$ norm of the solution $u(t)$ at the right of the slow soliton $R_2(t)$,
and the functional $g(t)$ corresponds locally to the norm adapted to each soliton. In particular, we have $g(t)\leq \|\eta(t)\|_{H_c^1}$.

We expand $u(t)=R_1(t)+R_2(t)+\eta(t)$ in the three quantities $\int u^2(t)$, $\mathcal{I}(t)$ and $E(u(t))$. 

\begin{lemma}[Expansion of energy type quantities]\label{dd} For all $t\in [0,t^*]$,
\begin{equation}\label{dd1}
\left| \int u^2(t)  - \left(c_1^{2q}(t)+c_2^{2q}(t)\right)\int Q^2 -\int\eta^2(t) \right| \le K e^{-\frac {\sqrt{c}}4 t} \exp(-2 c^{-r}),
\end{equation}
\begin{equation}\label{dd2}
\left|\mathcal{I}(t) - c_1^{2q}(t) \int Q^2 -\int\eta^2(t) \psi(x-m(t)) \right| \le K e^{-\frac {\sqrt{c}}{32} t} \exp(-2 c^{-r}), 
\end{equation}
\begin{equation}\label{dd3}\begin{split}
&\left|E(u(t))-\left\{ E(R_1)+E(R_2)+{ \frac 12} \int \eta_x^2(t)-p\left(R_1^{p-1}(t)+R_2^{p-1}(t)\right)\eta^2(t)  \right\}\right|
\\Ê& 
\quad \le K D_0 (\alpha +  \exp(-{\textstyle\frac 12}c^{-r})) g(t)+ K e^{-\frac {\sqrt{c}}4 t} \exp(-2 c^{-r}) ,
\end{split}\end{equation}
\begin{equation}\label{dd5}
\left| E(R_j(t))-E(R_j(0))
+\frac {c_j(0)} 2 \left[c_j^{2q}(t)-c_j^{2q}(0)\right] \int Q^2
\right| \le K c_j^{2q+1}(0) \left(\frac {c_j^{2q}(t)}{c_j^{2q}(0)}-1\right)^2.
\end{equation}
\end{lemma} 
 Lemma \ref{dd} is proved in Appendix A.2. In the rest of this section, we assume
$\alpha$ and $c$ small enough so that 
\begin{equation}\label{centieme}
\|\eta(t)\|_{H^1}+|c_1(t)-1|+|\tfrac {c_2(t)}{c}-1|+ D_0 (\alpha + \exp(-{\textstyle \frac 12 r})) \leq \frac 1{100}.
\end{equation}

\medskip

We next obtain a contradiction  from the following lemma.

\begin{lemma}\label{APRIORICLAIM} There exists $D_0>0$ such that
for $\alpha,c>0$ small enough, independent of $t^*$,\begin{equation}\label{aprioriun}
\sup_{t\in [0,t^*]} \|u(t)-Q(.-\rho_1(t))-Q_c(.-\rho_2(t))\|_{H^1_c}\le \tfrac 12 {D_0} (\alpha c^{q+\frac 12}+\exp(-c^{-r})).
\end{equation}
\end{lemma}
If $t^*<+\infty$, then Lemma  \ref{APRIORICLAIM} and the continuity in $H^1$ of $u(t)$ contradict the definition of $t^*$.
Therefore, we only have to prove Lemma \ref{APRIORICLAIM}.

\medskip

\noindent\emph{Proof of Lemma \ref{APRIORICLAIM}.}
\noindent\emph{Step 1.}  Monotonicity result on $\mathcal{I}(t)$.

\begin{claim}[Almost monotonicity property of $\mathcal{I}$ ]\label{MONOTONIE}
For $\alpha$ and $c$ small enough, 
\begin{equation}\label{monot}
\forall t\in [0,t^*],\quad 
\mathcal{I}(t)-\mathcal{I}(0)\le K \exp(-  c^{-\frac 12 -r}).\end{equation}
\end{claim}
\noindent\emph{Proof of Claim \ref{MONOTONIE}.} 
By standard calculations, we have
\[\begin{split}
 {\mathcal{I}}'(t)&= - 3 \int u_x^2 \psi'(x{-}m(t)) - m'(t)\int u^2 \psi'(x{-}m(t)) +\int u^2 \psi'''(x{-}m(t)) \\ &
\quad + \frac{2p}{p+1} \int u^{p+1}\psi'(x{-}m(t))\\
& \leq  -\frac 3{20} \int u^2 \psi'(x{-}m(t)) + \frac{2p}{p+1} \int u^{p+1}\psi'(x{-}m(t))
\end{split}\]
(we have used $\psi'>0$, \eqref{surphiff} and $m'(t)\ge 2/5$ by  \eqref{treize}).

In the nonlinear term $\int u^{p+1}\psi'(x{-}m(t))$, we expand $u(t)=R_1(t)+R_2(t)+\eta(t)$. We obtain
\begin{equation*}\begin{split}
\int u^{p+1}\psi'(x{-}m(t)) & \leq K' \int u^2 \left(R_1^{p-1}+R_2^{p-1}+|\eta|^{p-1}\right) \psi'(x{-}m(t)) \\
&\leq 
 K' \int u^2 R_1^{p-1} \psi'(x{-}m(t))+  \frac 1{10} \int u^2 \psi'(x{-}m(t)) 
\end{split}\end{equation*}
for $\alpha$ and $c$ small enough, since $\|\eta\|_{L^\infty}^{p-1}\le K \|\eta\|_{H^1}^{p-1}\leq K c$ and $0<R_2^{p-1}\leq K c$.
Moveover, by calculations similar to the ones of Claim \ref{prems}  and $\|u\|_{L^\infty}\le K$, we have
\[
\int u^2  R_1^{p-1}  \psi'(x{-}m(t))  \le   K   e^{-\frac t {32} } \exp(-  c^{-\frac 12 -r}).
\]
Thus, for all $t'\in [0,t^*]$,
$
{\mathcal{I}}'(t') \le K   c e^{-\frac {t'} {32} } \exp(-  c^{-\frac 12 -r}).
$
Let $t\in [0,t^*]$. By integration on $[0,t]$, we obtain 
$
{\mathcal{I}}(t)-{\mathcal{I}}(0)\le K \exp(-  c^{-\frac 12-r}).
$\hfill  $\Box$

 \medskip
 
\noindent\emph{Step 2.} Estimates on the scaling  parameters.
Let
\begin{equation*}
\Delta_j(t)=\frac {c_j^{2q}(t)}{c_j^{2q}(0)}-1.
\end{equation*}

\begin{claim}\label{DELTA} For all $t\in [0,t^*]$,
\begin{equation}\label{delta1}
|\Delta_1(t)|+c^{2q+1}|\Delta_2(t)| \le K (g(t)+g(0)+ \exp(-2 c^{-r})).
\end{equation}\end{claim}

\noindent\emph{Proof  of Claim \ref{DELTA}.}
Since there are only two solitons, the proof follows only from the $L^2$ norm and the energy conservation, i.e.
\eqref{dd1}, \eqref{dd3} and \eqref{dd5}. (When there are more than three solitons,  the use of quantities such as $\mathcal{I}(t)$ is also needed, see \cite{MMT}.)
Let $t\in [0, t^*]$.

From \eqref{dd1} taken at time $0$ and $t$, and $\int u^2(t)=\int u^2(0)$,   
\begin{equation}\label{dd1bis}\begin{split}
|c_1^{2q}(0) \Delta_1(t)+c_2^{2q}(0)\Delta_2(t)| & \le K \left(\int \eta^2(t) +  \int \eta^2(0) + \exp(-2 c^{-r})\right) \\
& \leq \frac K c (g(t)+g(0) + \exp(-2 c^{-r})).
\end{split}\end{equation}
From \eqref{dd3}, $E(u(t))=E(u(0))$ and
\begin{equation*}
 	\left| \int \eta_x^2(t)-p\left(R_1^{p-1}(t)+R_2^{p-1}(t)\right)\eta^2(t) \right|
	\leq K g(t),
\end{equation*}
we have
$$|E(R_1(t))-E(R_1(0))+E(R_2(t))-E(R_2(0))|\leq K (g(t)+g(0) + \exp(-2 c^{-r})).$$
Then, from \eqref{dd5}, we obtain
\begin{equation}\label{dd3bbis}
|c_1^{2q+1}(0) \Delta_1(t)+c_2^{2q+1}(0)\Delta_2(t)| \le K (g(t)+g(0)) + \Delta_1^2(t) + c^{2q+1} \Delta_2^2(t)
+ K\exp(-2c^{-r}).
\end{equation}
Multiplying \eqref{dd1bis} by $c_2(0)$ and combining with \eqref{dd3bbis}, from \eqref{centieme},  we obtain
\[
c_1^{2q}(0)(c_1(0)-c_2(0)) |\Delta_1(t)| \le K (g(t)+g(0)) +\Delta_1^2(t) + c^{2q+1} \Delta_2^2(t) + K\exp(-2c^{-r}).
\]
By \eqref{centieme},  we obtain
\[
|\Delta_1(t)| \le K (g(t)+g(0)) +\Delta_1^2(t) + c^{2q+1} \Delta_2^2(t) + K\exp(-2c^{-r}).
\]
Using this estimate in \eqref{dd3bbis}, we obtain similarly
\[
c^{2q+1}|\Delta_2(t)|  \le K (g(t)+g(0)) +\Delta_1^2(t) + c^{2q+1} \Delta_2^2(t) + K\exp(-2c^{-r}).
\]
Therefore, for $\Delta_1(t), \Delta_2(t)$ small enough (by \eqref{onze}), we obtain
\[
|\Delta_1(t)|+c^{2q+1}|\Delta_2(t)| \le K (g(t)+g(0))+K \exp(-2 c^{-r}).
\]
\quad\hfill$\Box$

\medbreak

\noindent\emph{Step 3.} Main argument of the proof of stability. 

For $t\in [0,t^*]$, as in \cite{MMT} we set
\[
\mathcal{F}(u(t))=E(u(t))+\frac {c_2(0)}2 \int u^2(t) + \frac {c_1(0)-c_2(0)}2 \mathcal{I}(t).
\]
The functional $\mathcal{F}$ coincides in a neighborhood of $R_1$ (respectively, $R_2$)
 with the functional introduced by Weinstein in \cite{We1} to prove the  stability of $R_1$ (resp., $R_2$).

We claim the following result on the quadratic part (in $\eta$)  of $\mathcal{F}(u(t))$.

\begin{claim}\label{quadra} Let
\begin{equation*} 
H(t)
=\frac 12  \int \eta_x^2(t)+ \left[c_2(0)  + (c_1(0){-}c_2(0)) \psi(x{-}m(t))\right]\eta^2(t) -p\left(R_1^{p-1}(t)+R_2^{p-1}(t)\right)\eta^2(t).\end{equation*}
There exists $\lambda_0>0$ independent of $c$ such that, for all $t\in [0,t^*]$,
\begin{equation}\label{qq1}
 {\lambda_0} g(t)
\leq 
H(t) \le   \frac 1 {\lambda_0} g(t).
\end{equation}
\end{claim}
 See Appendix A.3 for the proof of Claim \ref{quadra}.
 
\medskip

On the one  hand, using \eqref{dd1}-\eqref{dd5}, we obtain the following estimate
\begin{equation}\label{MMjams}\begin{split}
|\mathcal{F}(u(t))-\mathcal{F}(u(0)) - (H(t)-H(0))| &\le K D_0 (\alpha + \exp(-{\textstyle \frac 12} c^{-r})) (g(t)+g(0)) 
\\Ê&\quad + K \Delta_1^2(t) + K c^{2q +1} \Delta_2^2(t) + K \exp(-2 c^{-r}).
\end{split}\end{equation}
By $H(0)\leq K g(0)$, we obtain 
\begin{equation}\label{MMjamsbb}\begin{split}
|\mathcal{F}(u(t))-\mathcal{F}(u(0)) - H(t)| &\le K g(0)
+ K D_0 (\alpha + \exp(-{\textstyle \frac 12} c^{-r})) g(t)
\\Ê&\quad + K \Delta_1^2(t) + K c^{2q +1} \Delta_2^2(t) + K \exp(-2 c^{-r}).
\end{split}\end{equation}

On the other hand, by conservation of $E(u(t))$ and $\int u^2(t)$, and by the monotonicity of $\mathcal{I}(t)$
(see \eqref{monot}), we have
\begin{equation}
\mathcal{F}(u(t))-\mathcal{F}(u(0)) \le K \exp(-  c^{-\frac 12-r}),
\end{equation}
and thus by Claims \ref{quadra} and \ref{DELTA},
we obtain
\begin{equation}\label{MMarchives}\begin{split}
g(t) & \leq \frac 1{\lambda_0} H(t) \leq \frac 1{\lambda_0} \left( \mathcal{F}(u(t))-\mathcal{F}(u(0))
 + |\mathcal{F}(u(t))-\mathcal{F}(u(0)) -H(t)|\right) \\  
&
 \leq K g(0)+ K \big[D_0 (\alpha + \exp(-{\textstyle \frac 12} c^{-r})) + \Delta_1(t) + \Delta_2(t)\big] g(t) + K \exp(-2 c^{-r}).
\end{split}\end{equation}
For $\alpha$ and $c$ small enough, $g(t)$, $\Delta_1(t)$, $\Delta_2(t)$ and $D_0 (\alpha + \exp(-{\textstyle \frac 12} c^{-r}))$ are small, and
from \eqref{quatorze},  we obtain the following.

\begin{claim}\label{pourstab}
\begin{equation}\label{MMnonlinearity}
\|\eta(t)\|_{H^1_c}^2\le g(t)\le K g(0) + K \exp(- 2 c^{-r})
\leq K \alpha^2 c^{2q+1}  + K \exp(- 2 c^{-r}).
\end{equation}
\begin{equation}\label{mm5}
|\Delta_1(t)|+c^{2q+1}|\Delta_2(t)| \leq K (g(t)+g(0) + \exp(- 2 c^{-r}))
\le K \alpha^2 c^{2q+1}  + K \exp(- 2 c^{-r}).
\end{equation}
\begin{equation}\label{f4-54}
|c_1(t)-1|+c^{q+\frac 12} \left| \frac {c_2(t)}{c} -1 \right|
\leq K \alpha c^{q+\frac 12} + K \exp(- 2 c^{-r}).
\end{equation}
\end{claim}
Note that \eqref{mm5} follows from Claim \ref{DELTA} and \eqref{MMnonlinearity}.

\medskip
 
Now, we go back to $u(t)$ to prove \eqref{huit}.
From  direct calculations (recall that $\|Q_c\|_{H^1_c}\sim K c^{q+\frac 12}$) 
and \eqref{f4-54}, we have
\[
\|R_1(t)-Q(.-\rho_1(t))\|_{H^1_c} \le K |c_1(t)-1| \le K \alpha c^{q+\frac 12} + K \exp(- 2 c^{-r}),
\]
\[
\|R_2(t)-Q_{c}(.-\rho_2(t))\|_{H^1_c} \le K c^{q+\frac 12}\left|\frac {c_2(t)} c-1\right| 
\le K \alpha c^{q+\frac 12}  + K \exp(- 2 c^{-r}).
\]
Therefore,  we obtain
\begin{eqnarray*}&&
\|u(t)-Q(.-\rho_1(t))-Q_c(.-\rho_2(t))\|_{H^1_c}
\\&& \quad
\le \|u(t)-R_1(t)-R_2(t)\|_{H^1_c}+\|R_1(t)-Q(.-\rho_1(t))\|_{H^1_c}+\|R_2(t)-Q_c(.-\rho_2(t))\|_{H^1_c}
\\&& \quad
\le \|\eta(t)\|_{H^1_c} + K\alpha c^{q+\frac 12}+
 K \exp(- c^{-r})
 \le K_1 \left(  \alpha c^{q+\frac 12}  +\exp(- c^{-r})\right),\end{eqnarray*}
 where $K_1$ is independent of $\alpha$, $c$ and $D_0$.
 
 Choose now 
 \[
 D_0=4K_1,\]
  and then choose $\alpha>0$, $c>0$ small enough, so  that all the previous estimates hold. Then,
 for all $t\in [0,t^*]$, we have
 \[
 \|u(t)-Q(.-\rho_1(t))-Q_c(.-\rho_2(t))\|_{H^1_c} \le \frac 14{D_0} \left(  \alpha c^{q+\frac 12}  +\exp(- c^{-r})\right).
\qquad \Box \]

\section{Refined asymptotics for the $2$-soliton structure}
We claim the following.

\begin{proposition}[Asymptotic stability]\label{ASYMPTOTICVRAI}
There exist $K >0$, $\alpha_0>0$, $\overline c_0>0$  such that for any $0<c<\overline c_0$, $0<\alpha<\alpha_0$ the following if true.

Let $u(t)$ be an $H^1$ solution of \eqref{fkdv} such that
for $\frac 12 T_c \leq X_0 \leq \frac 32 T_c$,
\begin{equation}
\label{D100}
\|u(0)- Q-Q_c(.+X_0 )\|_{H^1}\le \alpha c^{q+\frac 12},
\end{equation}
so that Proposition \ref{ASYMPTOTIC} applies with   $\rho_1(t)$, $\rho_2(t)$. Then
\begin{enumerate}
\item Convergence of $u(t)$.
There exist $c_1^+, c_2^+ >0$  such that
\begin{equation}\label{neuf}
\lim_{t\rightarrow +\infty}\|u(t)-Q_{c_1^+}(x-\rho_1(t))-Q_{c_2^+}(x-\rho_2(t))\|
_{H^1(x> ct/10)}=0.
\end{equation}
\begin{equation}\label{sept}
|c_1^+-1|\leq   K  \alpha c^{q+\frac 12 }+ K \exp(-c^{- \frac 1{400}}),\quad
 \left|\frac {c_2^+}{c} - 1\right|\le K  \alpha  + K \exp(-c^{- \frac 1{400}}),
\end{equation}
\item  
Assume further that 
$f(u)=u^p$, for $p=2,3$ or $4$, and
 $\int_{x>0} x^2 \, u^2(0,x) dx <K_0$. Then, there exist  $x_1^+$ and 
$x_2^+$ such that
\begin{equation}
\lim_{t\to +\infty} \rho_1(t)-c_1^+ t=x_1^+, \quad
\lim_{t\to +\infty} \rho_2(t)-c_2^+ t=x_2^+.
\end{equation}
In the case $p=4$, if in addition, for some $\kappa>0$,
\begin{equation}\label{decayfort}
\alpha< \kappa c^{\frac 13} \quad \text{and}\quad \int_{x>\frac {11} {12} |\ln c|} x^2 \, u^2(0,x) dx <  \kappa c^{\frac 54}
\end{equation}
then
\begin{equation}
|x_1^+-\rho_1(0)|\leq K c^{\frac 58}, \quad |x_2^+-\rho_2(0)|\leq K c^{\frac 1{12}}.
\end{equation}
\end{enumerate}
\end{proposition}

\noindent\emph{Remark.} To obtain the convergence of the translation parameters, one has to add
an extra assumption on the initial data such as \eqref{decayfort}. Indeed, in the energy space,
one can construct an explicit example where convergence does not hold (see \cite{MMnonlinearity}).
Condition \eqref{decayfort} is enough for our purposes and could  be relaxed, and adapted
for the cases $p=2,3$.

\medskip

In what follows, we concentrate on the case $f(u)=u^p$ for $p=2,3$ or $4$.
The proof of the asymptotic stability (part 1 of Proposition \ref{ASYMPTOTICVRAI})
in the case of a general nonlinearity $f(u)$ follows from 
\cite{MMT} and \cite{MMas1}. Note that estimate \eqref{sept} is a direct consequence of
\eqref{f4-54}.

In the proof of Proposition \ref{ASYMPTOTICVRAI},
we need   another proof of the asymptotic stability 
for $f(u)=u^p$,   for $p=2,3$ or $4$, which is derived from the direct arguments of
\cite{MMnonlinearity}. The interest of this direct approach is to obtain an estimate on the 
convergence (see Lemma \ref{INTETA}), which is fundamental in proving the convergence of the translation parameters.
For a general nonlinearity, this kind of property is open.

\medskip

\noindent\emph{Proof of  Proposition \ref{ASYMPTOTICVRAI}.} 

1. The argument presented now is very similar to \cite{MMnonlinearity}, proof of Theorem 1, Step 3.
We keep the notation of the proof of Proposition \ref{ASYMPTOTIC}, in particular, the 
decomposition of $u(t)$ introduced in Claim \ref{LEMMEB1} and the conclusion of Claim \ref{pourstab}.
Now, we prove that $c_1(t)$ and $c_2(t)$ converge as $t\to +\infty$, and
that $\eta(t)$ converges to $0$ in $H^1(x>ct/10)$ as $t\to +\infty$.

We first control $\eta(t)$ around the solitons.

\begin{lemma}[Asymptotic stability locally in space]\label{INTETA}
Let
\begin{equation*}\begin{split}
& g_1(t)=  \int  (\eta_x^2 +  \eta^2)(t,x) e^{-\frac 1 4\,|x-\rho_1(t)|}dx, \quad
g_2(t)=  \int  (\eta_x^2  + c \, \eta^2)(t,x) e^{-\frac {\sqrt{c}} 4\,|x-\rho_2(t)|}dx.
\end{split}\end{equation*}
Then,
$$
\int_0^{+\infty} \left(g_1(t)+c^{\frac 32} g_2(t)\right) dt \leq K \alpha^2 c^{2q+1} + K \exp(-2c^{-r}).
$$
\end{lemma}

\noindent\emph{Proof of Lemma \ref{INTETA}.}
The proof is based on a localized Viriel type estimate.
Consider $\Phi:\mathbb{R}\to \mathbb{R}$ be an even smooth function such that
\begin{equation*}\begin{split}
& \Phi(x)=1\hbox{ on }[0,1];\  \Phi(x)=e^{-x}\hbox{ on }[2,+\infty);\  e^{-x}\le \Phi(x)\le 3e^{-x},
 \ \Phi'\le 0 \hbox{ on }\mathbb{R}^+,\\
 &
\Psi(x)=\int_{0}^x \Phi(y) dy, \ L_0 = \int_0^{+\infty} \Phi, \ \Psi_A(x)=A\Psi\left(\frac xA\right)
\text{ for $A>0$,}\\
&
\Psi_j(t,x)= \Psi_{A}\left(\sqrt{c_j(t)}(x-\rho_j(t))\right),\ 
\Theta_1(t,x)=\Psi_1(t,x)+ L_0 A,\ 
\Theta_2(t,x)=\Psi_2(t,x)- L_0 A.
\end{split}
\end{equation*}
Remark that $\Theta_1>0$, $\Theta_1'>0$, $\lim_{-\infty} \Theta_1=0$ and
$\Theta_2<0$, $\Theta_2'<0$, $\lim_{+\infty} \Theta_2=0$. Let
\begin{equation*}
\mathcal{K}_1(t)=L_0 A\int R_1^2(t)+\int \Theta_1(t)\eta^2(t),\quad
\mathcal{K}_2(t)=-L_0 A\int R_2^2(t)+\int \Theta_2(t)\eta^2(t).
\end{equation*}

\begin{claim}[Viriel estimate]\label{VIRIEL}There exist  $A\geq 5$, $K>0$,  $\alpha_0>0$ and $\overline c_0>0$ such that for $0<\alpha<\alpha_0$, $0<c<\overline c_0$  and   for all $t\in [0,+\infty)$,
\begin{equation}\label{viriel}
 \sqrt{c_j(t)} g_j(t) \leq K \left(-\frac {d}{dt} \mathcal{K}_j(t)
+ e^{-\frac 18\sqrt{c}(t+T_c)}\right).
\end{equation}
\end{claim}
See the proof in Appendix B.1. Note that this result is very similar to Lemma 2 in 
\cite{MMnonlinearity}. 
In \cite{MMnonlinearity},  the identity was  established in the case of one soliton. Here, to treat the two soliton situation, we have to use an additional term in
$\mathcal{K}_1(t)$, $\mathcal{K}_2(t)$.

\medskip

From now on, we fix $A\geq 5$ such  that Claim \ref{VIRIEL} holds. Then
$
|\Theta_1(x)|+|\Theta_2(x)|\leq K=K(A).
$
From \eqref{mm5} in Claim \ref{pourstab}, 
and $\int R_j^2= c_j^{2q}(t) \int Q^2$, we have, for $t\geq 0$,
\begin{equation*}
|\mathcal{K}_1(t)-\mathcal{K}_1( 0)|\le K (|c_1^{2q}(t)-c_1^{2q}( 0)|+\|\eta(t)\|_{L^2}^2
+\|\eta( 0)\|_{L^2}^2)\le K (\alpha^2  {c^{2q+1}}  +    \exp(-2 c^{-r})),
\end{equation*}
\begin{equation*}
|\mathcal{K}_2(t)-\mathcal{K}_2( 0)|\le K(|c_2^{2q}(t)-c_2^{2q}( 0)|+\|\eta(t)\|_{L^2}^2
+\|\eta( 0)\|_{L^2}^2)\le K (\alpha^2 {c^{2q}} +    \exp(-2 c^{-r})).
\end{equation*}
Therefore, integrating (\ref{viriel}) on $[0,+\infty)$ and using   $T_c=c^{-\frac 12-
\frac 1{100}}$, we obtain ($r=\frac 1{400}$),
\begin{equation}\label{conv1}
\int_{0}^{+\infty} g_1(t) dt  \le 
K (\alpha^2 c^{2q+1} +   \exp(-2 c^{-r})),
\end{equation}
\begin{equation}\label{conv2}
\int_{0}^{+\infty}  \sqrt{c} \, g_2(t) dt  \le 
K (\alpha^2 c^{2q} +   \exp(-2 c^{-r})).
\end{equation}
Thus Lemma \ref{INTETA} is proved.

\medskip

Now, we control  the scaling parameters.
Estimate (\ref{geom2}) and Lemma \ref{INTETA} imply that $c_1(t)$ and $c_2(t)$ have limits as $t\rightarrow +\infty$, which we denote respectively by $c_1^+$ and $c_2^+$.
By the stability result \eqref{f4-54},
\begin{equation*}
|c_1^+-1|\leq   K  (\alpha c^{q+\frac 12 }+   \exp(-c^{- r})),\quad
 \left|\frac {c_2^+}{c} - 1\right|\le K ( \alpha  + \exp(-c^{- r})),
\end{equation*}

Now, we extend the convergence of $\eta$ to $0$ in a large region in space, following
the proof of Theorem 1 in \cite{MMnonlinearity}.
We give a sketch the proof (see \cite{MMnonlinearity}, proof of Theorem 1, Step 3, for more details).

From (\ref{conv1}), (\ref{conv2}) (Viriel argument), there exists a sequence $(t_n)$ with $t_n\in [n,n+1)$ such that $g_1(t_n)+g_2(t_n)\to 0$ as $n\to +\infty$. Using $g_j'(t)\leq K_j g_j(t)$ (by a direct computation using \eqref{eqeta}, \eqref{geom1}, \eqref{geom2}), we obtain
$\lim_{t\to +\infty} (g_1(t)+g_2(t))=0$.

The rest of the proof is based only on monotonicity arguments on $u(t)$ and $u_x(t)$ such as in Claim \ref{MONOTONIE}, applied on different regions. Set
\begin{equation}\label{pourI}
\mathcal{I}_{\sigma,y_0}(t)=\int (u_x^2+u^2)(t,x) \psi(\sqrt{\sigma}(x-\tfrac \sigma 4 -y_0)) dx. 
\end{equation}
First, for $x_0>0$, let $y_0=\rho_1(t_0) -\frac \sigma 4 t_0 +x_0$ and $\sigma=1$. Using
$\mathcal{I}_{\sigma,y_0}(t_0)-\mathcal{I}_{\sigma,y_0}(0)$, we obtain
\begin{equation*}
\limsup_{t\rightarrow +\infty} \int_{x>\rho_1(t)+x_0} (\eta_x^2+\eta^2)(t,x)dx\le K e^{-\frac{x_0}{16}}.
\end{equation*} 
Next, for $y_0=\rho_2(t_0)-\frac \sigma 4 t_0+x_0$, $\sigma=c$, using
 $\mathcal{I}_{\sigma,y_0}(t_0)-\mathcal{I}_{\sigma,y_0}(\bar t_0)$, with
 $\frac c4 \bar t_0 +y_0=\rho_1(\bar t_0)-x_0$ (note that $\bar t_0\geq \frac c{10} t_0$ for $t_0$ large),
 we deduce
\begin{equation*}
\limsup_{t\rightarrow +\infty} \int_{x>\rho_2(t)+x_0}(\eta_x^2 +c \eta^2)(t,x)dx\le K e^{-\sqrt{c}\frac{x_0}{16}}.
\end{equation*} 
Finally,  by applying another monotonicity argument for
$y_0=0$ and $\sigma = \frac 25c$, we obtain
\begin{equation}\label{vers0}
\lim_{t\rightarrow +\infty} \int_{x>\frac c{10} t}(\eta_x^2 +  \eta^2)(t,x)dx=0.
\end{equation} 
Estimate \eqref{neuf} follows.

\medskip

2. Now, we prove the second part of Proposition \ref{ASYMPTOTICVRAI}.
Assume that 
\begin{equation}\label{decay}
\int_{x>0} x^{2} \, u^2(0,x) dx <K_0 
\quad \text{and}\quad \gamma_0(\alpha,c)= \int_{x>|\ln (\alpha c^{q+\frac 12})|} x^2 \, u^2(0,x) dx<1.\end{equation}
Note that
\begin{equation}\label{ast}
\rho_j(t)- \rho_j(0) -c_j^+ t = 
\int^t_0 (c_j(s)-c_j^+) ds + \int_0^t (\rho'_j(s)-c_j(s)) ds. 
\end{equation}
To prove that $\rho_j(t)-c_j^+t$ has a limit as
$t\to +\infty$, we will study separately the existence of  limits as $t\to +\infty$ of the two integrals above.

\medskip

\noindent\emph{2a. Preliminary : Monotonicity results on $\eta(t)$.}
We introduce monotonicity results on $\eta(t)$ (and not on $u(t)$ as before) that are refinement of  Claim \ref{MONOTONIE}.

We define, for $j=1,$ $2$,
\begin{equation}\label{defME}
\begin{split}
& \mathcal{M}_j(t)=\int \eta^2 \psi_j ,\\
&\mathcal{E}_j(t)=\int  \left[\frac 12 \eta_x^2
-\frac 1{p{+}1}\left( (R_1{+}R_2{+}\eta)^{p+1} {-}(p{+}1) R_1^p \eta {-} (p{+}1) R_2^p \eta {-}(R_1{+}R_2)^{p+1}\right)\right]
\psi_j,
\end{split}
\end{equation}
where $\psi_1(x)=\psi(\tilde x)$,
$\tilde x= x-\rho_1(t)+\frac 12 (t-t_0) $
and $\psi_2(x)=\psi(\sqrt{c}\tilde x_c)$, $\tilde x_c= x-\rho_2(t)+ \frac c2 (t-t_0) $, for $t_0\geq 0 $
and  $\psi(x)$ is defined by  \eqref{surphiff}.

We claim the following monotonicity results (see the proof in Appendix B.2).
\begin{claim}\label{MONOETA}
 For all $t\geq t_0\geq 0$,
\begin{equation*}
\begin{split}
	& \frac d{dt}\left( c_1^{2q}(t) \int Q^2+ \mathcal{M}_1(t) \right) 
	   \leq K e^{-\frac 1{16}(t-t_0)} g_1(t) + Ke^{-\frac 1 {32} \sqrt{c} (t+T_c)},\\
	& \frac d{dt}\left(-\frac {2q}{2q+1} c_1^{2q+1}(t) \int Q^2+ 2 \mathcal{E}_1(t) 
	+\frac 1{100}  \left(  c_1^{2q}(t) \int Q^2+ \mathcal{M}_1(t)  \right)\right) 
	 \\ &  \qquad \leq K e^{-\frac 1{16}(t-t_0)} g_1(t)+ Ke^{-\frac 1 {32} \sqrt{c} (t+T_c)}.\\
	   & \frac d{dt}\left( \left(c_1^{2q}(t)+c_2^{2q}(t)\right) \int Q^2+ \mathcal{M}_2(t) \right) 
	   \leq K e^{-\frac {c\sqrt{c}}{16}    (t-t_0)} \sqrt{c}\, g_2(t) + Ke^{-\frac 1 {32} \sqrt{c} (t+T_c)},\\
	& \frac d{dt}\left(-\frac {2q}{2q{+}1} \left(c_1^{2q+1}(t){+}c_2^{2q+1}(t)\right) \int Q^2+ 2 \mathcal{E}_2(t) 
	+\frac c{100} \left( \left(c_1^{2q}(t){+}c_2^{2q}(t)\right) \int Q^2+ \mathcal{M}_2(t) \right) \right) \\ &
	  \qquad  \leq K e^{-\frac {c\sqrt{c}}{16}   (t-t_0)  } c^{\frac 32} g_2(t)+ Ke^{-\frac 1 {32} \sqrt{c} (t+T_c)}.
\end{split}
\end{equation*}
\end{claim}

\noindent\emph{Remark.} 
The improvement of Claim \ref{MONOETA} with respect to the monotonicity on $u(t)$ in Claim \ref{MONOTONIE}, or Lemma 3 in \cite{MMnonlinearity} is that the upper bound can be integrated twice in time.

Claim \ref{MONOETA} is   one example of monotonicity  result on localized energy type quantities
on $\eta(t)$. In Appendix B.2, we prove a slightly more general version of Claim \ref{MONOETA}, where
$\tilde x= x-\rho_1(t)+\frac \sigma 2 (t-t_0)-x_0$, where $0<\sigma\leq \frac 12$ and
$x_0\in \mathbb{R}$, and we claim other monotonicity results to be used in this paper.

\medbreak

Let 
\begin{equation*}
\tilde g_1(t)=\int \big(\eta_x^2{+}\eta^2\big)(t,x) \psi(x{-}\rho_1(t)) dx,
\quad
\tilde g_2(t)=\int \big(\eta_x^2{+} c \eta^2 \big)(t,x) \psi(\sqrt{c}(x{-}\rho_2(t))) dx.
\end{equation*}
Note that there exists $K >0$ such that $\psi(x)\geq \frac 1 K e^{-\frac {|x|}4}$ on $\mathbb{R}$. Thus,  we have
$g_j(t)\leq K \tilde g_j(t)$.

We have the following consequence of Claim \ref{MONOETA}.
\begin{claim}[Control on the scaling parameters]\label{SURLESC}
	For all $t\geq t_0$,
\begin{equation}\label{surlescun}
\begin{split}
 & |c_1(t)-c_1(t_0)|   \leq K (\tilde g_1(t) + \tilde g_1(t_0)) + K \int_{t_0}^t e^{-\frac 1{16} (t'-t_0)} g_1(t') dt' + Ke^{-\frac 1 {32} \sqrt{c} (t_0+T_c)},\\
 &|c_2^{2q+1}(t)-c_2^{2q+1}(t_0)|  \leq 
2 |c_1(t)-c_1(t_0)|
+K (\tilde g_2(t)+\tilde g_2(t_0))
 \\ & + K c^{\frac 32} \int_{t_0}^t e^{-\frac {c \sqrt{c}}{16} (t'-t_0) } g_2(t') dt'
+ Ke^{-\frac 1 {32} \sqrt{c} (t_0+T_c)}.
\end{split}
\end{equation}
\end{claim}

\noindent\emph{Proof of Claim \ref{SURLESC}.}
Integrating the conclusion of Claim \ref{MONOETA} between $t_0$ and $t$, we obtain
\begin{equation*}
\begin{split}
&   c_1^{2q}(t)- c_1^{2q}(t_0) \leq \frac 1{\int  Q^2}
(\mathcal{M}_1(t_0)-\mathcal{M}_1(t))
  +K\int_{t_0}^t e^{-\frac 1{16} (t'-t_0)} g_1(t')dt'
+	  Ke^{-\frac 1 {32} \sqrt{c} (t_0+T_c)},\\
&   \frac {2q}{2q+1} \left(c_1^{2q+1}(t)- c_1^{2q+1}(t_0) \right)-\frac 1{100}\left(c_1^{2q}(t)- c_1^{2q}(t_0)\right)
 \geq   \frac 1{\int  Q^2} \left(2\mathcal{E}_1(t)-2\mathcal{E}_1(t_0)\right)
 \\ & \quad +\frac 1{100} \frac 1{\int  Q^2} \left(\mathcal{M}_1(t)-\mathcal{M}_1(t_0)  \right)
  -K\int_{t_0}^t e^{-\frac 1{16} (t'-t_0)} g_1(t')dt'
  -  Ke^{-\frac 1 {32} \sqrt{c} (t_0+T_c)}.
\end{split}
\end{equation*}
Since  $c_1\sim 1$ by \eqref{centieme}, $\mathcal{M}_1(t)\geq 0$, 
$\mathcal{E}_1(t)\geq - K g_1(t) \geq -K' \, \tilde g _1(t)$, 
and $\mathcal{E}_1(t_0)\leq K \, \tilde g_1(t_0)$, $\mathcal{M}_1(t_0)\leq K \, \tilde g_1(t_0)$,
we obtain the first estimate of \eqref{surlescun}. 
The estimate on $|c_2^{2q+1}(t)-c_2^{2q+1}(t_0)|$ is obtained in the same way
using $\mathcal{M}_2(t)$ and $\mathcal{E}_2(t)$.
\quad $\Box$

\medskip

We claim the following lemma.

\begin{lemma}\label{OdeT}
Assume that \eqref{decay} holds.  
Then,
$$\int_0^{+\infty} \tilde g_1(t) dt  
\leq K  ((\alpha c^{q+\frac 12})^{\frac 74} +  \gamma_0(\alpha,c)+
\alpha^{- \frac 14}\exp(-\tfrac 32 c^{-r}) ) ),
$$
$$
\int_0^{+\infty} \tilde g_2(t) dt \leq K (\alpha^{\frac 74} c^{\frac 74 q-\frac 18}  + \alpha^2 c^{\frac 32q -\frac 12} + \tfrac 1 c \gamma_0(\alpha,c)+  \alpha^{- \frac 14}\exp(-\tfrac 32 c^{-r})),$$
$$\lim_{t\to +\infty} t \, \int (\eta_x^2(t,x)+\eta^2(t,x))\psi(x-\tfrac{c}{10} t) dx=0.$$
\end{lemma}
For the proof  see Appendix B.3. Note that the proof is based only on Lemma \ref{INTETA} and monotonicity
arguments such as Claim \ref{MONOETA}. We follow the same steps as in the proof 
of \eqref{vers0}, using   quantities  $\mathcal{M}_j(t),$ $\mathcal{E}_j(t)$ instead
of $\mathcal{I}_{\sigma, y_0}$ on the same lines.
The  proof of the monotonicity is  the same as  the one of Claim \ref{MONOETA}.

\medskip

\noindent\emph{2b. Estimate on $\rho_j(t)-c_j^+$.}

\begin{lemma}[Estimate on $c_j(t)-c_j^+$]\label{params}
Assume that \eqref{decay} holds.  
Then,
\begin{equation*}\begin{split}
& \int_0^{+\infty} |c_1(t)-c_1^+| dt \leq K  ((\alpha c^{q+\frac 12})^{\frac 74} +  \gamma_0(\alpha,c) +
\alpha^{- \frac 14} \exp(- c^{-r})) ,\\
& \int_0^{+\infty} |c_2(t)-c_2^+| dt  \leq K (\alpha^{\frac 74} c^{-\frac 14 q-\frac 18}  + \alpha^2 c^{-\frac 12 q -\frac 12} + c^{-2q-1} \gamma_0(\alpha,c)+ \alpha^{ - \frac 14} \exp(- c^{-r})).
\end{split}\end{equation*}
\end{lemma}
\noindent\emph{Proof of Lemma \ref{params}.}
Let $t\to +\infty$ in Claim \ref{SURLESC}, since
$c_j(t)\to c_j^+$ and $\tilde g_j(t)\to 0$, we obtain
\begin{equation*}
|c_1^+ - c_1(t_0)|\leq K \, \tilde g_1(t_0)   + K \int_{t_0}^{+\infty}
e^{-\frac 1{16} (t'-t_0) } g_1(t') dt' + Ke^{-\frac 1 {32} \sqrt{c} (t_0+T_c)}.
\end{equation*}
\begin{equation*}
c^{2q} |c_2^+-c_2(t_0)|\leq K |c_1^+ - c_1(t_0)|+ 
K \, \tilde g_2(t_0)   + K c^{\frac 32} \int_{t_0}^{+\infty}
e^{-\frac {c  \sqrt{c} } {16} (t'-t_0)} g_2(t') dt' + Ke^{-\frac 1 {32} \sqrt{c} (t_0+T_c)}.
\end{equation*}
Thus, integrating on $[0,+\infty)$ and using Fubini Theorem, 
since $g_1(t)\leq K \tilde g_1(t)$, we obtain
\begin{equation*}
\int_0^{+\infty} |c_1(t)-c_1^+| dt \leq K \, \int_0^{+\infty} \tilde g_1(t) dt   
+  K\alpha^{- \frac 14}\exp(-2 c^{-r}).
\end{equation*}
Similarly,
\begin{equation*}\begin{split}
& c^{2q} \int_0^{+\infty} |c_2(t)-c_2^+| dt  \leq K \, \int_0^{+\infty} (\tilde g_1(t)+ \tilde g_2(t)) dt +  K\alpha^{- \frac 14} \exp(-2 c^{-r}).
\end{split}\end{equation*}
Thus Lemma \ref{params} is a consequence of Lemma \ref{OdeT}.\quad $\Box$

\begin{lemma}[Estimate on $\rho_j'(t)-c_j(t)$]\label{semiconv}
Assume that \eqref{decay} holds.

For $j=1,2$,  $\int_0^{+\infty} (\rho_j'(t)-c_j(t)) dt$ is defined
and
\begin{equation*}\begin{split}
 	\biggl|	\int_0^{+\infty} (\rho_1'(t)-c_1(t)) dt \biggr|& \leq 
K((\alpha c^{q+\frac 12})^{\frac 78}+\gamma_0^{\frac 12}(\alpha,c)+  \exp(-c^{-r})),
\\ 	\biggl|	\int_0^{+\infty} (\rho_2'(t)-c_2(t)) dt \biggr|   &\leq K(\alpha^{\frac 78} c^{-\frac 14  q+\frac 12}
 +\alpha^2 c^{-\frac 12}
+ \alpha^3 c^{-\frac 12-\frac 14 q}\\ & + (c^{-q+\frac 14} +\alpha^2 c^{-q -\frac 34})
 \gamma_0^{\frac 12}(\alpha,c) +  \exp(-\tfrac 14c^{-r})).
\end{split}\end{equation*}
\end{lemma}
Lemma \ref{semiconv} is proved in Appendix B.4. Note that it makes use of the following
functional
\begin{equation}\label{defJ}
	J_j(t)=c_j^{-2q}(t)\int \eta(t,x) \left(\int_{-\infty}^x \tilde R_j(t,x')dx'\right) dx,
\end{equation}
where $\tilde R_j(t,x)=\tilde Q_{c_j(t)}(x-\rho_j(t))$, which is an $L^1$-type quantity, already
introduced in \cite{MMjams}.
For
$p=3$, another argument can be used.
From  Claim \ref{GEOM}, \eqref{geom4},
$$|\rho_j'(t)-c_j(t)|\leq K g_j(t) + K e^{-\frac 18 \sqrt{c} (t+T_c)},$$
where $g_j(t)$ is defined in \eqref{defgj}. By \eqref{conv1}-\eqref{conv2},
we obtain in this case:
$$
(p=3)\quad 
\int_0^{+\infty} |\rho_1'(t)-c_1(t)| + c^{\frac 32}  |\rho_2'(t)-c_2(t)| dt \leq 
K \alpha^2 c^{2q+1} + K \exp(-2c^{-r}).
$$
Such an integrability property cannot be proved from \eqref{geom1} for $p=2, 4$,
since we do not know whether or not
$\int_0^{+\infty} \sqrt{g_j(t)}dt <+\infty.$

\medskip

From \eqref{ast}, Lemmas \ref{params} and \ref{semiconv}, it follows that 
for $j=1,\ 2$, $$\rho_j(t)-c_j^+ t \to x_j^+\quad \text{ as $t\to +\infty$,}$$
where 
\begin{equation}\label{ast2}\begin{split}
 |x_1^+ -\rho_1(0)| &\leq K((\alpha c^{q+\frac 12})^{\frac 78}+\gamma_0^{\frac 12}(\alpha,c)+ 
\alpha^{- \frac 14}\exp(-c^{-r})),\\
|x_2^+ -\rho_2(0)|&\leq  
K(\alpha^{\frac 78} c^{-\frac 14  q+\frac 12}+\alpha^{\frac 74} c^{-\frac 14 q-\frac 18}  + \alpha^2 c^{-\frac 12 q -\frac 12} +  c^{-2q-1} \gamma_0(\alpha,c)
\\&\  
 + (c^{-q+\frac 14} +\alpha^2 c^{-q -\frac 34})
 \gamma_0^{\frac 12}(\alpha,c) +\alpha^{-\frac 14} \exp(-\tfrac 14c^{-r})).
\end{split}\end{equation}

\medskip

\noindent\emph{3. Control of the shifts under assumption \eqref{decayfort} for $p=4$.}
Now, we assume, for some $\kappa>1$,
\begin{equation}\label{decayfort2}
\gamma_0=\gamma_0(\kappa c^{\frac 13}, c) \leq \int_{x>\frac {11}{12}|\ln c|} x^2 \, u^2(0,x) dx < K c^{\frac 54}.
\end{equation}
We apply the previous estimate with $\bar \alpha= \kappa c^{\frac 13}$. Then,  
from \eqref{ast2}, using $q=\frac 1{12}$ in this case, we get
\begin{equation*}\begin{split}
  |x_1^+ -\rho_1(0)|\leq K c^{\frac 78(\frac 13+\frac 1{12}+\frac 12)}+
K c^{\frac 58} + K c^{-\frac 1{12}} \exp(-\tfrac 12c^{-r})) \leq K c^{\frac 58},
\quad |x_2^+ -\rho_2(0)|\leq  Kc^{\frac 1{12}},
\end{split}\end{equation*}
where the worst term in $|x_2^+ -\rho_2(0)|$ is $c^{-2q-1} \gamma_0(\alpha,c))\leq K c^{\frac 1 {12}}$.

\appendix

\section{Appendix}

\subsection{Proof of Claim \ref{LEMMEB1}}
We first  state a preliminary result. Recall $T_c=c^{-\frac 12 (1+\frac 1{100})}$.
For $\alpha, c>0$, we define
\begin{equation}\label{33}\begin{split}
\mathcal{U}(\alpha,c)=& \Big\{u\in H^1(\mathbb{R}); \text{ $\exists r_1,r_2 \in \mathbb{R}$ $|$
$|r_1-r_2|>{ \tfrac 12} T_c$ } \\& \ \text{and }
\Big\|u- Q(.-r_1)-Q_{c}(.-r_{2})\Big\|_{H^1}\le \alpha c^{q}\Big\}.\end{split}
\end{equation}

\begin{lemma}[Existence of modulation parameters]\label{MODULATION}
 There exists $\overline c_{0}>0$, $\alpha_0>0$, $K>0$ and unique
$C^1$ functions $(c_1,c_2,\rho_1,\rho_2):{\mathcal U}(\alpha_0,\overline c_{0})
\rightarrow (0,+\infty)^2\times \mathbb{R}^2 ,$ such that if $u\in {\mathcal U}(\alpha_0,\overline c_{0})$, and 
\begin{equation}\label{32}
\eta(x)=u(x)- Q_{c_1}(x-\rho_1)-Q_{c_2}(x-\rho_2),
\end{equation}
then, for $j=1,2$,
\begin{equation}\label{orthoA}
\int Q_{c_j}(x-\rho_j) \,\eta(x)dx=
\int Q_{c_j}'(x-\rho_j) \,\eta(x)dx=0.
\end{equation}
Moreover, if $u\in
{\mathcal U}(\alpha,c)$, with $0<\alpha<\alpha_0$, $0<c<\overline c_{0}$,
then
\begin{equation}\label{smallnessA}
\|\eta\|_{H^1}+|c_1-1|\le K \alpha c^{q } ,\quad \left|\frac {c_2}{c}-1\right|\le K \alpha,
\quad |\rho_1-\rho_2|>{\textstyle\frac 14} T_c.
\end{equation}
\end{lemma}

\noindent\emph{Proof.} 
Let $\alpha_0$, $\overline c_{0}>0$ to be chosen later.
Let $0<\alpha<\alpha_0$, $0<c<\overline c_{0}$. 
First, let $r_1,r_2\in\mathbb{R}$ be such that $|r_1-r_2|\ge \frac 12 T_c$ and consider
\begin{equation}
\mathcal{V}(\alpha)=\mathcal{V}_{c,r_1,r_2}(\alpha)=
\bigl\{u\in H^1(\mathbb{R}), \|u-Q(.-r_1)-Q_c(.-r_2)\|_{H^1}\le 2 \alpha c^{q}\big\}.
\end{equation}
Let $V(\alpha)=[0,\alpha]^2\times [-\alpha,\alpha]^2\times \mathcal{V}(\alpha)\subset \mathbb{R}^4\times H^1(\mathbb{R})$,
\begin{equation}
M_0=(0,0,0,0,Q(.-r_1)+Q_c(.-r_2)),\text{ and any } 
M=(\lambda_1,\lambda_2,y_1,y_2,u)\in V(\alpha).
\end{equation}
Let also
\begin{equation*}
Q_1(x)=Q_{1+\lambda_1 c^{q }}(x-r_1-c^{q } y_1 ),\quad
Q_2(x)=Q_{c(1+\lambda_2 )}(x-r_2-c^{-\frac 12} y_2),
\end{equation*}
\begin{equation*}
w(M)=c^{-q}(u-Q_1-Q_2),
\end{equation*}
\begin{eqnarray*}
&& \nu_1(M)=\int w(M) Q_1,\quad \nu_2(M)=c^{-q} \int w(M)Q_2,\\
&& \mu_1(M)=\int w(M) Q_{1}',\quad \mu_2(M)=c^{-(q+\frac 12)}\int w(M)Q_{2}'.
\end{eqnarray*}
 For any $M\in V(\alpha)$, 
since $\int Q_c^2=c^{2q} \int Q^2$ and $\int (Q_{c}')^2=c^{2(q+\frac 12)}\int (Q')^2$, we have
\begin{eqnarray}\label{borne}
|\nu_1(M)|+
|\nu_2(M)|+|\mu_1(M)|+|\mu_2(M)]\le K\alpha. 
\end{eqnarray}
\begin{claim}\label{MODUL} For any $M\in V(\alpha)$, for any $j,k=1,2$, $j\neq k$,
\begin{eqnarray*} &&
\left|\frac {\partial \nu_j}{\partial \lambda_j}(M)+\frac {5{-}p}{4(p{-}1)} \int Q^2\right|+\left|\frac {\partial \mu_j}{\partial y_j}(M)-\int (Q')^2\right|
+\left|\frac {\partial \nu_j}{\partial y_j}(M)\right|
+\left|\frac {\partial \mu_j}{\partial \lambda_j}(M)\right|\le K\alpha,
\\    &&
\left|\frac {\partial \nu_j}{\partial y_k}(M)\right|
+\left|\frac {\partial \mu_j}{\partial \lambda_k}(M)\right|
+\left|\frac {\partial \nu_j}{\partial \lambda_k}(M)\right|
+\left|\frac {\partial \mu_j}{\partial y_k}(M)\right|
\le K\exp(- c^{ -r}).
\end{eqnarray*}\end{claim}

\noindent\emph{Proof of Claim \ref{MODUL}.}
The claim follows from elementary calculations similar to the ones in \cite{MMT}, Appendix A.
We give the proof of some of these estimates. 
First, note that
\[
\frac {\partial \nu_1}{\partial \lambda_1}(M)=- c^{-q} \int \frac {\partial Q_1}{\partial \lambda_1} Q_1 + \int w(M) \frac {\partial Q_1}{\partial \lambda_1}.
\]
Moreover, by $Q_{c_0}(x)=c_0^{\frac 1{p-1}} Q(\sqrt{c_0} x)$, we have (recall 
$\tilde Q_{c_0}=\frac 2{p-1} Q_{c_0} + x Q_{c_0}'$)
\[  \frac {\partial Q_{c_0}}{\partial c_0}
= \frac {\tilde Q_{c_0}}{2c_0}
\quad\text{and  thus}\quad
  \int  \frac {\partial Q_{c_0}}{\partial c_0}  Q_{c_0}=\frac {5-p}{4(p-1)} c_0^{2q-1} \int Q_{c_0}^2.
\]
We have
\[
\frac {\partial Q_1}{\partial \lambda_1}=c^q \frac 1{2(1+\lambda_1 c^q)} \tilde Q_{1+\lambda_1 c^q}.
\]
Thus,
\[
\frac {\partial \nu_1}{\partial \lambda_1}(M)=- \frac {5-p}{4(p-1)} (1+\lambda_1 c^q)^{2q-1} \int Q^2 + \int w(M) \frac {\partial Q_1}{\partial \lambda_1},
\]
and by $\|w\|_{L^2}\le K \alpha$ and $|\lambda_1|\le \alpha$, we obtain
\[
\left|\frac {\partial \nu_1}{\partial \lambda_1}(M) + \frac {5{-}p}{4(p{-}1)} \int Q^2 \right|\le K\alpha.
\]
Similarly, 
\[
\frac {\partial \nu_2}{\partial \lambda_2}(M)=- c^{-2q} \int \frac {\partial Q_2}{\partial \lambda_2} Q_2 + c^{-q}\int w(M) \frac {\partial Q_2}{\partial \lambda_2}.
\]
We have
\[
\frac {\partial Q_2}{\partial \lambda_2}=  \frac 1{2(1+\lambda_2)} \tilde Q_{c(1+\lambda_2)}
\quad \text{and so}\quad
c^{-2q} \int \frac {\partial Q_2}{\partial \lambda_2} Q_2
=(1+\lambda_2)^{2q-1} \frac {5{-}p}{4(p{-}1)} \int Q^2.
\]
From 
\[
\left|c^{-q} \int  w(M) \frac {\partial Q_2}{\partial \lambda_2} \right|\le \|w\|_{L^2} c^{-q} \left\| \frac {\partial Q_2}{\partial \lambda_2}\right\|_{L^2} \le K \|w\|_{L^2}^2,
\]
and $|\lambda_2|\le K\alpha$, we obtain
\[
\left|\frac {\partial \nu_2}{\partial \lambda_2}(M) + \frac {5{-}p}{4(p{-}1)} \int Q^2 \right|\le K\alpha.
\]

Third, we have
\[
\frac {\partial \mu_2}{\partial y_2}(M)=- c^{-2q-\frac 12} \int \frac {\partial Q_2}{\partial y_2} Q_2' + c^{-q}\int w(M) \frac {\partial Q_2'}{\partial y_2}.
\]
But
\[
\frac {\partial Q_2}{\partial y_2}=-c^{-\frac 12} Q_2',\quad  \frac {\partial Q_2'}{\partial y_2}=- c^{-\frac 12} Q_2'',
\]
and so
\[
- c^{-2q-\frac 12} \int \frac {\partial Q_2}{\partial y_2} Q_2'=c^{-2(q+\frac 12)} \int (Q_2')^2
=(1+\lambda_2)^{2(q+\frac 12)} \int (Q')^2.
\]
Moveover,
\[
c^{-q} \left| \int w \frac {\partial Q_2'}{\partial y_2}\right| \le K \alpha.
\]

Fourth, we have
\[
\frac {\partial \nu_2}{\partial y_2}(M)=c^{-2q} \int \left(-\frac {\partial Q_2}{\partial y_2}\right) Q_2
+c^{-q} \int w  \frac {\partial Q_2}{\partial y_2}.
\]
The first term in the right hand side is 0  by parity, the second term is controlled by $K\alpha$.

Finally, we check a different term:
\[
\frac {\partial \nu_1}{\partial y_2}(M)=   c^{-q-\frac 12} \int Q_2' Q_1,
\]
since $\frac {\partial Q_1}{\partial y_2}=0$. By $\left|\int Q_2' Q_1\right|\le K   \exp(-2 c^{-r})$ (see proof of Claim \ref{prems} for similar estimates), we obtain the desired estimate. \hfill$\Box$

\bigbreak

By Claim \ref{MODUL} and (\ref{borne}),  and $(\nu_1,\nu_2,\mu_1,\mu_2)(M_0)=(0,0,0,0)$, we
apply the implicit fonction theorem to 
$(\nu_1,\nu_2,\mu_1,\mu_2)$: there exists $\overline c_0>0$, $\alpha_0>0$
(chosen independent of $c$, $r_1$, $r_2$) such that, if $0<c<\overline c_0$, $0<\alpha<\alpha_0$,
for all $u\in \mathcal{V}(\alpha)$, there exists $\lambda_1(u),\lambda_2(u),y_1(u),y_2(u)$ unique
such that if $M(u)=(\lambda_1(u),\lambda_2(u),y_1(u),y_2(u),u)$, then
\begin{equation}
\nu_1(M(u))=\nu_2(M(u))=\mu_1(M(u))=\mu_2(M(u))=0.
\end{equation}
Moreover,
\begin{equation}\label{gggg}
|\lambda_1(u)|+|\lambda_2(u)|+|y_1(u)|+|y_2(u)|\le K\alpha.
\end{equation}
Now, we set 
\begin{equation*}
c_1(u)=1+\lambda_1(u)c^{q},   c_2(u)=c(1+\lambda_2(u)), 
y_1(u)=r_1+y_1(u) c^{q},   y_2(u)=r_2+y_2(u) c^{-\frac 12}.
\end{equation*}

For any $u\in {\mathcal U}(\alpha,c)$, there exist $r_1,r_2$ satisfying $|r_1-r_2|\ge \frac 12 T_c$
and $\|u-Q(.-r_1)-Q(.-r_2)\|_{H^1}\le 2\alpha c^q$.
Thus $c_1(u), c_2(u), \rho_1(u), \rho_2(u)$ are defined as before for such $u$.
Uniqueness and regularity are consequences of the  implicit function theorem.
Note finally that (\ref{gggg}) implies
$|c_1-1|\le K\alpha c^{q},$ $|\frac {c_2} c- 1|\le K\alpha$ 
and $|x_1-x_2|\ge \frac 14 T_c$.
\hfill$\Box$

\bigskip

\noindent\emph{Proof of  Claim \ref{LEMMEB1}.}
For $t=0$, using assumption \eqref{48bon}, we apply Lemma \ref{MODULATION}  to $u(0)$ with $\alpha c^{\frac 12}$ instead of $\alpha$. We find $\rho_1(0)$, $\rho_2(0)$, $c_1(0)$ and $c_2(0)$ such that \eqref{dix} and \eqref{quatorze} hold.

For $t\in (0,t^*]$, by the definition of $t^*$, $\|u(t)\|_{H^1}\leq c^{-\frac 12}\|u(t)\|_{H^1_c}\leq D_0 (\alpha +c^{-\frac 12} \exp(- c^{-r}))$. We apply Lemma \ref{MODULATION} to $u(t)$ where $\alpha$ is replaced by $D_0 (\alpha +c^{-\frac 12} \exp(- c^{-r}))$ which is small for $\alpha$ small depending on $D_0$. 
We obtain directly \eqref{dix}--\eqref{onze}. The estimates on $\rho_1'(t)$ and $\rho_2'(t)$ follow from the equation of $\eta(t)$ written in Claim \ref{LEMMEB1}, see  
Claim \ref{GEOM} below. Since $\rho_1'(t)-\rho_2'(t)\geq \frac 12$ and $\rho_1(0)-\rho_2(0)\geq \frac 12 T_c$, we obtain \eqref{treize}.
\hfill$\Box$

\medbreak

\subsection{Proof of Lemma \ref{dd}}
First, we recall well-known identities related to $Q_c$ for $f(u)=u^p$.

\begin{claim}[Identities for any $p>1$]\label{LemmaA1}
    \begin{equation*}
        \int Q^{p+1} = \frac{2(p+1)}{p+3} \int Q^2, \quad       \int (Q')^2 = \frac{p-1}{p+3} \int Q^2.
    \end{equation*}
    \begin{equation*}
    		\int Q_c^2=c^{2q} \int Q^2,\quad E(Q_c)=c^{2q+1}E(Q)=-\frac {5-p}{2(p+3)} c^{2q+1} \int Q^2.
    \end{equation*}
\end{claim}

\noindent\emph{Proof of Lemma \ref{LemmaA1}.} These are well-known calculations. We have
$Q^p=Q-Q''$ and $\frac 2{p+1} Q^{p+1} = Q^2 -(Q')^2$. Thus, by integration:
$$
\int Q^{p+1} =\int Q^2 +\int (Q')^2,\quad
\frac 2{p+1} \int Q^{p+1} = \int Q^2 - \int (Q')^2.
$$
Therefore, $\int Q^{p+1} = \frac {2(p+1)}{p+3} \int Q^2$ and
  $\int (Q')^2=\int Q^{p+1}-\int Q^2= \frac {p-1}{p+3} \int Q^2$. Moreover,
$E(Q)= \frac 12 \int (Q')^2 -\frac 1{p+1} \int Q^{p+1}= \frac {p-5}{2(p+3)} \int Q^2$.

Since $Q_c(y)=c^{\frac 1{p-1}} Q(\sqrt{c} y)$ and $q= \frac 1{p-1} -\frac 14$, 
we have
$$\int Q_c^2(y) dy=c^{\frac 2{p-1}} \int Q^2(\sqrt{c} y) dy = c^{2q} \int Q^2.$$
Similary, $\int (Q_c')^2=c^{2q+1} \int (Q')^2$
and $\int Q_c^{p+1}= c^{2q+1} \int Q^{p+1}$, and so $E(Q_c)=c^{2q+1} E(Q)$.
\hfill $\Box$

\medskip

Then, we claim the following estimates (recall $r=\frac 1{400}$).
\begin{claim}\label{prems} For all $t\in [0,t^*]$,
\begin{equation}\label{cc0}
\forall x\in \mathbb{R},\quad R_1(t,x) \le K \psi(x-m(t));
\end{equation}
\begin{equation}\label{cc1}
0\le \int R_1(t,x)R_2(t,x)dx \le K e^{-\frac {\sqrt{c}}4 t} \exp(-2 c^{-r});
\end{equation}
\begin{equation}\label{cc2}
0\le \int R_1(t,x)(1-\psi(x-m(t)))dx \le K e^{-\frac t {32}} \exp(-c^{-\frac 12 -r});
\end{equation}
\begin{equation}\label{cc3}
0\le \int R_2(t,x)\psi(x-m(t)) dx \le K e^{-\frac {\sqrt{c}}{32} t} \exp(-2 c^{-r}).
\end{equation}
\end{claim}

\noindent\emph{Proof of Claim \ref{prems}.}
First, note that since $c_1(t)>\frac 12$ and $c_2(t)> \frac c 2$, we have
$$
0\leq R_1(t,x)\leq K e^{-\frac 12 |x-\rho_1(t)|},\quad
0\leq R_2(t,x)\leq K e^{-\frac c2 |x-\rho_2(t)|}.
$$

Proof of \eqref{cc0}. For $x>m(t)$, we have  $\psi(x-m(t))>\frac 12$ and so $R_1(t,x)\le K \psi(x-m(t))$.  For $x\le m(t)$, by the definition of $\psi(x)$ and $m(t)<\rho_1(t)$, we have
\[
R_1(t,x)\le K e^{\frac 12 (x-\rho_1(t))} \le K e^{\frac 12 (x-m(t)) + \frac 12(m(t)-\rho_1(t))} \le K \psi(x-m(t)).
\]

Proof of \eqref{cc1}. We have
\[
\begin{split}
 0 \le R_1(t,x) R_2(t,x) & \le K e^{-\frac 12 |x-\rho_1(t)|} e^{-\frac {\sqrt{c}}2 |x-\rho_2(t)|}\\&
 \le K e^{-\frac 12 |x-\rho_1(t)|} e^{\frac {\sqrt{c}}2 |x-\rho_1(t)|}e^{-\frac {\sqrt{c}}2 |\rho_1(t)-\rho_2(t)|}\\&
 \le K e^{-\frac 14 |x-\rho_1(t)|} e^{-\frac {\sqrt{c}}4 t} e^{-\frac 18 c^{-2r} }
  \le K e^{-\frac 14 |x-\rho_1(t)|} e^{-\frac {\sqrt{c}}4 t} \exp(-2 c^{-r}) .
\end{split}
\]
for $c$ small enough. Thus by integration in $x$, we obtain \eqref{cc1}.

Proof of \eqref{cc2}. For $x\ge m(t)+\frac t 8 + \frac 1 {32} T_c$, 
\[
1-\psi(x-m(t))\le K e^{-\frac 14 (x-m(t))} \le K e^{-\frac t{32}} \exp(-{\textstyle \frac 1 {64}} T_c),\]
and so $
R_1(t) (1-\psi(x-m(t)))\le K e^{-\frac 12 |x-\rho_1(t)|}   e^{-\frac t{32}} \exp(-{\textstyle\frac 1 {64}} T_c).$

For $x\le m(t)+\frac t 8 + \frac 1 {32}T_c\le \rho_1(t)-\frac t 8 - \frac 1 {32} T_c$,  \[
R_1(t)\le K e^{-\frac 12 |x-\rho_1(t)|}\le K e^{-\frac 14 |x-\rho_1(t)|} e^{-\frac t{32}} \exp(-{\textstyle\frac 1 {128}}T_c),
\]
and $0\le 1-\psi(x-m(t))\le 1$.
Thus, by integration in $x$ and for $c$ small, we obtain \eqref{cc2}.

Proof of \eqref{cc3}. For $x\le m(t)-\frac t 8 - \frac 1 {32}T_c$, 
\[
\psi(x-m(t))\le K e^{\frac 12 (x-m(t))} \le K e^{-\frac t{32}} \exp(-{\textstyle \frac 1 {64}} T_c),
\]
\[
R_2(t) \psi(x-m(t))\le K c^{q+\frac 14} e^{-\frac {\sqrt{c}}2 |x-\rho_2(t)|}   e^{-\frac t{32}} \exp(-{\textstyle\frac 1 {64}} T_c).
\]
For $x\ge m(t)-\frac t 8 - \frac 1 {32} T_c\ge \rho_2(t)+\frac t 8 + \frac 1 {32} T_c$, we have
\[
R_2(t)\le Kc^{q+\frac 14} e^{-\frac {\sqrt{c}} 2 |x-\rho_2(t)|}\le K c^{q+\frac 14} e^{-  \frac {\sqrt{c}}4 |x-\rho_2(t)|} e^{-\sqrt{c} \frac t{32}} \exp(-{\textstyle\frac 1 {64}} c^{-2r}),
\]
and $0\le \psi(x-m(t))\le 1$.
Thus, by integration in $x$ and for $c$ small, we obtain \eqref{cc3}.
\hfill$\Box$

\medbreak

\noindent\emph{Proof  of Lemma \ref{dd}.} Proof of \eqref{dd1}.
From $\int R_j \eta=0$, we have
\[
\int u^2(t)=\int R_1^2(t)+\int R_2^2(t)+\int \eta^2(t)+2 \int R_1(t)R_2(t).
\]
Thus \eqref{dd1} is now a consequence of \eqref{cc1} and $\int R_j^2(t) = c_j^{2q}(t) \int Q^2$ (Claim \ref{LemmaA1}).

Proof of \eqref{dd2}. Using $\int R_1 \eta=0$, we have
\[\begin{split}
\mathcal{I}(t)&= \int u^2(t,x) \psi(x{-}m(t)) dx \\&
=\int R_1^2(t) +\int R_1^2(t) (1{-}\psi(x{-}m(t))) + \int \left[ R_2^2(t) +2 R_1(t)R_2(t)\right]
\psi(x{-}m(t))\\
& \quad  -2 \int \eta(t) R_1(t) (1{-}\psi(x{-}m(t))) + 2 \int \eta(t) R_2(t)\psi(x{-}m(t))
+\int \eta^2(t) \psi(x{-}m(t)).
\end{split}\]
Since $\|\eta(t)\|_{L^\infty}\le 1$, we can use estimates \eqref{cc2} and \eqref{cc3} to obtain \eqref{dd2}.

Proof of \eqref{dd3}. 
First, we prove the following estimate :
\begin{equation}\label{dd3bis}\begin{split}
&\left|E(u(t))-\left\{ E(R_1)+E(R_2)+{\frac 12} \int \eta_x^2(t)-p\left(R_1^{p-1}(t)+R_2^{p-1}(t)\right)\eta^2(t)  \right\}\right|\\Ê& 
\le K\left\{\int \left(R_1^{p-2}+R_2^{p-2}\right) |\eta|^3 + \int |\eta|^{p+1} \right\} + K e^{-\frac {\sqrt{c}}4 t} \exp(-2 c^{-r}).
\end{split}\end{equation}
Let $R(t)=R_1(t)+R_2(t)$ (note that $R^k\leq K(R_1^k+R_2^k)$). By expanding $u(t)$, we have
$
E(u(t))  = \frac 12 \int (R+\eta)_x^2 -\frac 1{p{+}1} \int (R+\eta)^{p+1} ,
$
and thus
\[
\begin{split}&
\left|E(u(t))-\left\{
 E(R) - \int \left(R''+R^p\right) \eta + \frac 12 \int \eta_x^2 - p R^{p-1} \eta^2\right\}\right| 
\\ & \le K \int \left(R^{p-2} |\eta|^3 +|\eta|^{p+1}\right)
 \le K \int \left(R_1^{p-2}+R_2^{p-2}\right) |\eta|^3 + \int |\eta|^{p+1}. 
\end{split}\]
By estimates similar to \eqref{cc1} related to the decay of $Q$ and its derivatives,
\[
|E(R)-E(R_1)-E(R_2)|\le K e^{-\frac {\sqrt{c}}4 t} \exp(-2 c^{-r}).
\]
Since $R_j''+R_j^p=c_j(t)R_j$ and $\int R_j \eta=0$, we have by \eqref{cc1}:
\[
\left|\int \left(R''+R^p\right) \eta \right| \leq \int |R^p-R_1^p-R_2^p| | \eta |
\le K e^{-\frac {\sqrt{c}}4 t} \exp(-2 c^{-r}).
\]
Similarly, we have
\[
\left|\int \left(R^{p-1}-R_1^{p-1}-R_2^{p-1}\right) \eta^2  \right| \le K e^{-\frac {\sqrt{c}}4 t} \exp(-2 c^{-r}).
\]
Thus, we obtain \eqref{dd3bis}.

Now, we continue proving \eqref{dd3} by estimating  $\int \left(R_1^{p-2}+R_2^{p-2}\right) |\eta|^3 + \int |\eta|^{p+1}$.
Let
\begin{equation}\label{13bis}
\beta=D_0(\alpha + c^{-q-\frac 12} \exp(-c^{-r})) \quad \text{so that by \eqref{onze}}\quad 
\|\eta(t)\|_{H^1}\le K  \beta c^q.
\end{equation}
Note that $\beta\le D_0(\alpha +  \exp(-{\textstyle \frac 12}c^{-r}))$, for $c$ small enough.
 We have by the Gagliardo--Nirenberg inequality:
\[\begin{split}
\int |\eta|^{p+1}& \le K \left(\int \eta_x^2\right)^{\frac {p-1}4} \left(\int \eta^2\right)^{\frac {p+3}4} \le K \beta^{p-1} \int \eta_x^2 + K \beta^{-\frac {(p-1)^2}{5-p}} \left(\int \eta^2\right)^{\frac {p+3}{5-p}}.
\end{split}\]
Since $\frac {p+3}{5-p}=1+\frac 1{2q}$,  $\frac 1 q-\frac {(p-1)^2}{5-p}=p-1$ ($2q=\frac {5-p}{2(p-1)}$) and using \eqref{13bis}, we obtain
\[
\beta^{-\frac {(p-1)^2}{5-p}} \left(\int \eta^2\right)^{\frac {p+3}{5-p}} 
=\beta^{-\frac {(p-1)^2}{5-p}} \left(\int \eta^2\right)^{\frac 1{2q}} \int \eta^2
\le  \beta^{\frac 1 q-\frac {(p-1)^2}{5-p}} c \int \eta^2 \leq  K \beta^{p-1}  c \int \eta^2
\]
Thus,
\[
\int|\eta|^{p+1} \le K \beta^{p-1} \left[\int \eta_x^2 +   c \int \eta^2\right].
\]
In addition, from \eqref{13bis}, 
\[
\int R_2^{p-2} |\eta|^3 \le \beta \int R_2^{p-1} \eta^2 + \beta^{-(p-2)} \int |\eta|^{p+1}
\le   K \beta \left[\int \eta_x^2 +  c \int \eta^2\right],
\]
\[
\int R_1^{p-2} |\eta|^3 \le K \|\eta\|_{L^\infty} \int R_1 \eta^2  \le K  \beta \int \eta^2 \psi(x-m(t)) .
\]
Gathering these estimates and \eqref{dd3bis}, we obtain \eqref{dd3} (note that $p\geq 2$).

Proof of \eqref{dd5}. By Claim \ref{LemmaA1}, we have
\[
E(R_j(t))-E(R_j(0))=
-\frac {5{-}p}{2(p{+}3)} \left( c^{2q+1}(t) - c^{2q+1}(0)\right) \int Q^2 .
\]
But since $\frac {2q+1}{2q}= \frac {p{+}3}{5{-}p}$,  
\[ \left| c_j^{2q+1}(t) - c_j^{2q+1}(0)- \frac {p{+}3}{5{-}p} c_j(0)\left[c^{2q}(t) - c^{2q}(0)\right]\right| \le K c_j^{2q+1}(0) \left(\frac {c_j^{2q}(t)}{c_j^{2q}(0)}-1\right)^2.
 \]
Thus \eqref{dd5} follows.
\hfill$\Box$

\subsection{Proof of  Claim \ref{quadra}}
The proof is based on the  following well-known fact: 
\textit{There exists $\lambda_1>0$ such that
if $v\in H^1(\mathbb{R})$ satisfies $\int  Q v=\int xQ v=0$, then}
\begin{equation}\label{single}
\int v_x^2-p Q^{p-1} v^2 + v^2\ge \lambda_1 \|v\|_{H^1}^2.
\end{equation}
 
It is similar to \cite{MMT}, Proof of Lemma 4.
Set
\[
H_0(t)=\int \left(\eta_x^2+  [c_2(t)+c_1(t)\psi(x{-}m(t))]\eta^2 - p  (R_1^{p-1}+R_2^{p-1}) \eta^2\right).
\]
Note that  $H_0(t)$ and $H(t)$   are easily compared. Indeed, we have
\[
|H(t)-H_0(t)| \le \int \left[|c_2(t){-}c_2(0)| + |c_1(t) {-}c_1(0){-}c_2(0)|\psi(x{-}m(t)) \right] \eta^2 
\]
Let $\epsilon_0>0$. By  \eqref{onze}, for $\alpha$ and $c$ small enough, we have
\[
|H(t)-H_0(t)|  \le \epsilon_0 \int [c+ \psi(x{-}m(t))] \eta^2.
\]
Thus, it is sufficient to prove \eqref{qq1} for $H_0(t)$ for some $\lambda_0>0$ independent of $c$ and $\alpha$. 
 
First, we consider a function $\Phi\in \mathcal{C}^2(\mathbb{R})$, $\Phi(x)=\Phi(-x)$, $\Phi'\le 0$ on $\mathbb{R}^+$, with
\[
\text{$\Phi(x)=1$ on $[0,1]$; $\Phi(x)=e^{-x}$ on $[2,+\infty)$, $e^{-x}\le \Phi(x) \le 3 e^{-x}$ on $\mathbb{R}^+$.}
\]
Let $\Phi_B(x)=\Phi(\frac x B)$. We recall the following claim from \cite{MMT}, page 355 (and references therein):
\textit{There exists $B_0>0$ such that, for all $B>B_0$, if $v\in H^1(\mathbb{R})$ satisfies $\int  Q v=\int xQ v=0$, then}
\begin{equation}\label{singlel}
\int \Phi_B(x) \left(v_x^2-p Q^{p-1} v^2 + v^2 \right)dx\ge \frac {\lambda_1} 4 \int \Phi_B(x) (v_x^2+v^2)dx.
\end{equation}
This result is a localized version of \eqref{single}, and is easily proved by direct calculations.

By a scaling argument, i.e. changing $x$ into $x\sqrt{c}$ and using the definition of $Q_c$, we have : 
If $v\in H^1(\mathbb{R})$ satisfies $\int  Q_c v=\int xQ_c v=0$, then
\begin{equation}\label{singlelc}
\int \Phi_{\frac B{\sqrt{c}}}(x) \left(v_x^2-p Q_c^{p-1} v^2 + c v^2 \right)dx
\ge \frac {\lambda_1} 4 \int \Phi_{\frac B{\sqrt{c}}}(x) (v_x^2+c v^2)dx.
\end{equation}

Now, we consider $\eta$ as in the proof of Lemma \ref{APRIORICLAIM}, i.e. satisfying the orthogonality conditions $\int \eta R_j(t)= \int \eta xR_j(t)=0$ for $j=1,2$.
Let $\Phi_1(t,x)= \Phi_{B}(x{-}\rho_1(t)) $ and $\Phi_2(t,x)= \Phi_{\frac B{\sqrt{c}}}(x{-}\rho_2(t))$.
We have
\[\begin{split}
H_0(t)& = \int \Phi_2 \left(\eta_x^2-p R_2^{p-1} \eta^2 + c_2(t) \eta^2 \right)dx
 + \int \Phi_1 \left(\eta_x^2-p R_1^{p-1} \eta^2 + c_1(t) \eta^2 \right)dx
\\ & \quad -p \int R_2^{p-1} (1-\Phi_2 ) \eta^2-p \int R_1^{p-1} (1-\Phi_1 ) \eta^2
\\ & \quad + \int (1-\Phi_1-\Phi_2) \eta_x^2+ \left[(1-\Phi_2) c_2(t)+ (1-\Phi_1) c_1(t)\psi(x{-}m(t))\right] \eta^2. 
\end{split}\]
Let $\lambda_2=\min(\frac {\lambda_1} 4,\frac 12)$. Since $1-\Phi_1-\Phi_2\ge 0$, $c_2(t)>\frac c2$ and $c_1(t)\ge \frac 12$, we have by \eqref{singlel} and \eqref{singlelc}:
\[\begin{split}
&  \int \Phi_2 \left(\eta_x^2-p R_2^{p-1} \eta^2 + c_2(t) \eta^2 \right)dx
+ \int \Phi_1 \left(\eta_x^2-p R_1^{p-1} \eta^2 + c_1(t) \eta^2 \right)dx
\\ &   + \int (1-\Phi_1-\Phi_2) \eta_x^2+ \left[(1-\Phi_2) c_2(t)+ (1-\Phi_1) c_1(t)\psi(x{-}m(t))\right] \eta^2
 \\ & \ge \lambda_2 \int (\eta_x^2 + [c+ \psi(x-m(t)) ]\eta^2).   
\end{split}\]
Finally, since $|R_2^{p-1}(t,x)|\le K c e^{-\frac {\sqrt{c}} 2 |x-\rho_2(t)|}$ and $\Phi_2(x)=0$ for $|x-\rho_2(t)|\le \frac B{\sqrt{c}}$, we have
\[
\int R^{p-1}_2 (1-\Phi_2) \eta^2 \le K c e^{-\frac {p-1} 2 B} \int \eta^2 \le \epsilon_0 c \int \eta^2,
\]
for $B$ large enough.
Similarly, for $B$ large enough, since $e^{-\frac 12 |x-\rho_1(t)|}\le K \psi(x-m(t))$, we have
\[
\int R^{p-1}_1 (1-\Phi_1) \eta^2 \le  K e^{-\frac {p-1} 2 B} \int \psi(x{-}m(t)) \eta^2\le \epsilon_0  \int \psi(x{-}m(t)) \eta^2.
\]
Therefore, for $B$ large enough, we obtain
\[
H_0(t)\ge \frac {\lambda_2} 2 \int \left(\eta_x^2 + [ c+ \psi(x{-}m(t))] \eta^2 \right).
\]
This completes the proof of Claim \ref{quadra}.
\hfill$\Box$

\section{Appendix}

\subsection{Proof of Claim \ref{VIRIEL} - Localized Viriel estimate}
By explicit calculations, $\eta(t)$ satisfies
\begin{equation}\label{eqeta}
\eta_t =(-\eta_{xx}-(R_1+R_2+\eta)^p+R_1^p+R_2^p)_x
-\frac {c_1'}{2c_1} \tilde R_1 -\frac {c_2'}{2c_2} \tilde R_2
 +(\rho_1'-c_1)R_{1x}+(\rho_2'-c_2)R_{2x},
\end{equation}
where $\tilde R_j(t,x)=\tilde Q_{c_j(t)}(x-\rho_j(t))$  and $\tilde Q_c(x)=\frac 2{p-1}Q_c + xQ_c'$.

\medskip

\textit{Step 1.}\quad Control of the geometrical parameters. 

\begin{claim}\label{GEOM} Let
\begin{equation}\label{defgj}\begin{split}
  g_1(t)=  \int  (\eta_x^2  +  \eta^2)(t,x) e^{-\frac 1 4\,|x-\rho_1(t)|}dx, \quad
  g_2(t)=  \int  (\eta_x^2 + c \, \eta^2)(t,x) e^{-\frac c 4\,|x-\rho_2(t)|}dx.
\end{split}\end{equation}
Then, for all $t\geq 0$,
\begin{equation}\label{geom1}
|(\rho_j'-c_j)c_j^{2q}|\le K c_j^{q+\frac 12}\sqrt{g_j(t)}
+ K g_j(t)+ K e^{-\frac 1{8}\sqrt{c}(t+T_c)},
\end{equation}
\begin{equation}\label{geom2}
|(c_j^{2q})'|\le K\sqrt{c_j} \, g_j(t)+ K e^{-\frac 1{8}\sqrt{c}(t+T_c)},
\end{equation}
\begin{equation}\label{geom4}
\left|\frac 12 (\rho_j'-c_j)c_j^{2q} \int Q^2 -(p-3) \int \eta R_j^p\right|
  \le 
K g_j(t)+K e^{-\frac 1{8}\sqrt{c}(t+T_c)}.
 \end{equation}
\end{claim}

\noindent\emph{Proof of Claim \ref{GEOM}.} 
Let $(j,k)=(1,2)$ or $(j,k)=(2,1)$. Since $\int Q \tilde Q=2q \int Q^2$ and $\int QQ_x=0$,
\begin{eqnarray}0&=&
\frac {d}{dt}\int \eta R_j =\int \eta_t R_j - \rho_j' \int \eta R_{jx} 
+\frac {c_j'}{2c_j} \int \eta \tilde R_j\nonumber \\ 
&=& -\int (-\eta_{xx}+c_j \eta-pR_j^{p-1} \eta) R_{jx} 
 -(\rho_j'-c_j) \int \eta R_{jx} +\frac {c_j'}{2c_j} \int \eta \tilde R_j \label{nul1}\\
&& +\int \left[(R_j+R_k+\eta)^p-R_j^p-R_k^p-pR_j^{p-1}\eta \right] R_{jx} 
\nonumber\\&&-q \frac {c_j'}{c_j} \int R_j^2 -\frac {c_k'}{2c_k} \int \tilde R_k R_j + (\rho_k'-c_k) \int R_{kx} R_j.\nonumber
\end{eqnarray}
First, we note that the first integral in (\ref{nul1}) is zero since
$\mathcal{L}$ is self-adjoint  and $\mathcal{L} Q'=0$.
Second, we have  $\sqrt{c_j}|R_j(t,x)|+|R_{jx}(t,x)|\le c^{\frac 1{p-1} +\frac 12} e^{-\sqrt{c_j}|x-\rho_j|}$,
and $\int R_j^2=c_j^{2q} \int Q^2$, $\int R_{jx}^2=c_j^{2q +\frac 12} \int Q_x^2$.
Finally,   since $\rho_1(t)-\rho_2(t)\geq \frac 12 (t+T_c)$, by the proof of Claim \ref{prems}, all the terms containing a product $R_j R_k$ or their derivatives, are controlled by $e^{-\frac 1{8}\sqrt{c}(t+T_c)}$.
Thus,
\begin{equation}\label{vraicp}
\begin{split}
& \biggl|\frac 12  (c_j^{2q})'(t) \int Q^2
  +(\rho_j'-c_j) \int \eta R_{jx} -\frac {c_j'}{2c_j} \int \eta \tilde R_j \\
& -\int \left[(R_j+\eta)^p-R_j^p-pR_j^{p-1}\eta \right] R_{jx} 
\biggr| \leq K e^{-\frac 1{8}\sqrt{c}(t+T_c)}.
\end{split}
\end{equation}

Next, we note that 
\begin{equation}\label{similaire}
|(R_j+\eta)^p-R_j^p-pR_j^{p-1}\eta |\le
K (|R_j|^{p-2} |\eta|^2 +|\eta|^p),
\end{equation}
and thus
\begin{eqnarray*} 
\left| (c_j^{2q})'  \right| &  \leq & 
K \left( c_j^{2q+\frac 12}|\rho_j'-c_j)|
+|(c_j^{2q})'| \right) \left[c_j^{-2q} \int \eta^2 e^{-\sqrt{c_j}|x-\rho_j|}\right]^{\frac 12}
\\ && + K \int  (|R_j|^{p-2} |\eta|^2 +|\eta|^p)|R_{jx}| + K e^{-\frac 1{8}\sqrt{c}(t+T_c)}.
\end{eqnarray*}
We claim 
\begin{equation}\label{encoreclaim}
 \int  (|R_j|^{p-2}  \eta ^2 +|\eta|^p)|R_{jx}|
  \leq  K  \sqrt{c_j} g_j(t) .
\end{equation}
Indeed, first,
\begin{equation*}\begin{split}
\int   |R_j|^{p-2}   |R_{jx}|\eta^2
\leq Kc_j^{\frac 32}\int \eta^2 e^{-\frac {\sqrt{c_j}}2(x-\rho_j)}
\leq K c_j^{\frac 12} g_j(t).\end{split}
\end{equation*}
Second (for $p=3$ and $4$),
 \begin{equation*}\begin{split}
 \int    |\eta|^p|R_{jx}|
  \leq c^{\frac 1{p-1}+\frac 12} \|\eta^2  \phi_{j}\|_{L^{\infty}}^{\frac {p-2}2}
   \int \eta^2 e^{-\frac {\sqrt{c_j}}2(x-\rho_j)}\end{split}
\end{equation*} 
where $\phi_{j}(x)=e^{-\frac {\sqrt{c_j}}2(x-\rho_j)}$. 
Since $ {\phi_{jx}}=-\frac{\sqrt{c_j}} 2  \phi_{j} $, we have
\begin{equation*}
\sup_{x\geq \rho_j(t)}|\eta^2(x) \phi_{j} (x)|^2 \leq \left[\int_{x\geq \rho_j(t)} \biggl|2\eta\eta_x \phi_{j}
+\eta^2 \frac {\phi_x}{2 \phi_{j}}\biggr|\right]^2\le
 K  \int \eta^2 \phi_{j}\left(\int_{x\geq \rho_j(t)} \eta_x^2\phi_{j} + c_j \eta^2 \phi_{j}\right).
\end{equation*}
Thus,
$\|\eta^2  \phi_{j}\|_{L^\infty(\mathbb{R})}^2 
\le 
 K \left( \int \eta^2 \phi_{j}\right) g_j(t).
$
We obtain 
 \begin{equation*}\begin{split}
 \int    |\eta|^p|R_{jx}|
 & \leq c^{\frac 1{p-1}+\frac 12} 
\left(   \int \eta^2 e^{-\frac {\sqrt{c_j}}2(x-\rho_j)}\right)^{\frac {p-2} 4 +1} g_j^{\frac {p-2}4}(t)
\\ & 
\leq K c^{\frac 1{p-1} - 1 +\frac p4} \left(   \int \eta^2 e^{-\frac {\sqrt{c_j}}2(x-\rho_j)}\right)^{\frac {p-2}2} g_j(t)
\leq K c^{-1+\frac p4 +\frac p {2(p-1)}} g_j(t)\leq K c ^{\frac 12 } g_j(t),\end{split}
\end{equation*} 
using
 $\|\eta\|_{H^1}\le K \alpha c^q$.
This proves (\ref{encoreclaim}).

Therefore, using again
 $\|\eta\|_{H^1}\le K \alpha c^q$, for $\alpha$ small,
\begin{equation}\label{peutservir1}
|(c_j^{2q})'|\le K|\rho_j'-c_j| c_j^{q+\frac 12}\left( \int \eta^2 e^{-{\sqrt{c_j}}(x-\rho_j)}\right)^{\frac 12}
+K \sqrt{c_j}\,g_j(t) + K e^{-\frac 1{8}\sqrt{c}(t+T_c)}.
\end{equation}

Using $\int \eta (x-\rho_j)R_j=0$ and $\int \eta R_j=0$, we have in a similar way
\begin{eqnarray*}0&=&
\frac d{dt}\int \eta (x-\rho_j)R_j=\int \eta_t (x-\rho_j)R_j-\rho_j' \int \eta  (x-\rho_j)R_{jx} +\frac{c_j'}{2c_j}\int \eta (x-\rho_j)  \tilde R_j \\
&=& -\int (-\eta_{xx}+c_j \eta-pR_j^{p-1} \eta)((x-\rho_j)R_{j})_x -(\rho_j'-c_j)\int \eta (x-\rho_j)R_{jx}\\
&&+ \frac {c_j'}{2c_j} \int \eta (x-\rho_j)\tilde R_j +\int \left[(R_j+R_k+\eta)^p-R_j^p-R_k^p-pR_j^{p-1}\eta \right]((x-\rho_j)R_j)_x\\
&& +(\rho_j'-c_j) \int R_{jx} (x-\rho_j)R_j -\frac{c_k'}{2c_k} \int \tilde R_k (x-\rho_j)R_j
+(\rho_k'-c_k)\int R_{kx}(x-\rho_j)R_j.
\end{eqnarray*}
From $\mathcal{L}((xQ)_x)=-2Q-(p-3)Q^p$ and (\ref{encoreclaim}), it follows that
\begin{eqnarray}&&
\left|\frac 12 (\rho_j'-c_j)c_j^{2q} \int Q^2 -(p-3) \int \eta R_j^p\right|\label{avantpeutservir2}\\&&
 \le Kg_j(t)+K \left(|\rho_j'-c_j| c_j^{2q}+c_j^{-\frac 12}|(c_j^{2q})'|\right) 
c_j^{-q}\left( \int \eta^2 e^{-\sqrt{c_j}|x-\rho_j|}\right)^{\frac 12}+ K e^{-\frac 1{8}\sqrt{c}(t+T_c)},
\nonumber\end{eqnarray}
which implies by (\ref{peutservir1}) and $\|\eta\|_{H^1}\le K \alpha c^q$,
\begin{equation}\label{peutservir2}
|(\rho_j'-c_j)c_j^{2q}|\le K c_j^{q+\frac 12}\sqrt{g_j(t)}
+ Kg_j(t)+ K e^{-\frac 1{8}\sqrt{c}(t+T_c)}.
\end{equation}
Now, inserting (\ref{peutservir2}) in (\ref{peutservir1}), we have the following estimate
\begin{equation}\label{finalcj}
|(c_j^{2q})'|\le K\sqrt{c_j}\, g_j(t)+ K e^{-\frac 1{8}\sqrt{c}(t+T_c)}.
\end{equation}
Finally, inserting (\ref{finalcj}) and (\ref{peutservir2}) in  (\ref{avantpeutservir2}),
we obtain \eqref{geom4}.
\hfill$\Box$

\medskip

\textit{Step 2.}\quad Viriel identity in $\eta$ and conclusion of the proof.
For $w\in H^1(\mathbb{R})$, define
\begin{equation*}
H_{\infty}(w,w)=\frac 12 \int \left(3w_x^2+ w^2- p \bigl(Q^{p-1}-(p-1) xQ_xQ^{p-2}\bigr)w^2\right),
\end{equation*}
\begin{equation*}
H_{\infty}^*(w,w)=H_{\infty}(w,w)-\frac{2(p-3)}{\int Q^2} \left(\int w xQ_x\right)\left(\int w Q^p\right).
\end{equation*}
Recall that the two quantities $H_{\infty}(w,w)$ and $H_{\infty}^*(w,w)$ were introduced in \cite{MMarchives}
for a Viriel type identity. Here, by analogy, we define, 
\begin{equation*}
H_j(\eta,\eta)=\frac 12 \int \left(3\eta_x^2\Phi_j+c_j \eta^2\Phi_j-p\bigl( R_j^{p-1}\Phi_j-(p-1)\Psi_j R_{jx}R_j^{p-2}\bigr)\eta^2\right),
\end{equation*}
\begin{equation*}
H_j^*(\eta,\eta)=H_j(\eta,\eta)-\frac{2(p-3)}{\int R_j^2} \left(\int \eta \Psi_j R_{jx}\right)\left(\int \eta R_j^p\right),
\end{equation*}
where $\Phi_j(t,x)=\frac \partial {\partial x} \Psi_j(t,x)=\sqrt{c_j(t)} \Phi\left(\frac {\sqrt{c_j(t)}}{A}(x-\rho_j(t))\right)$.

By  straightforward calculations,  we have
\begin{equation*}
 \frac {d}{dt}\mathcal{K}_j(t) =  \frac {d}{dt}\int \Theta_j(t)\eta^2(t)
+  \Theta_j(0)(c_j^{2q})'\int Q^2,
\end{equation*}
\begin{equation*}
\frac 12\frac {d}{dt}\int \Theta_j(t)\eta^2(t)=
\int \Theta_j \eta_t \eta- \frac {c_j}2 \int \Theta_{jx} \eta^2
-\frac 12(\rho_j'-c_j) \int \Theta_{jx} \eta^2+\frac {c_j'}{4c_j} \int (x-\rho_j) \Theta_{jx} \eta^2.
\end{equation*}
Therefore, by using the equation of $\eta$, we have, for $j\neq k$,
\begin{eqnarray*}
\frac 12 \frac {d}{dt}\mathcal{K}_j(t) &=&
-H_j^*(\eta,\eta)+\frac 12 \int \Phi_{jxx} \eta^2
-\frac 12 (\rho_j'-c_j) \int \Phi_j \eta^2+\frac {c_j'}{4c_j}\int (x-\rho_j)\Phi_j \eta^2\nonumber\\
&& +\int \left[(R_j+R_k+\eta)^p-R_j^p-R_k^p-pR_j^{p-1}\eta\right](\Theta_j \eta)_x\\
&& +\frac 12 \Theta_j(0)(c_j^{2q})'\int Q^2-\frac p2 \Theta_j(0)\int (R_j^{p-1})_x \eta^2+\Theta_j(0)(\rho_j'-c_j)\int \eta R_{jx}\\
&& -\frac {c_j'}{2c_j}\int \tilde R_j\Theta_j\eta-\frac {c_k'}{2c_k}\int \tilde R_k\Theta_j\eta+(\rho_k'-c_k)\int R_{kx}\Theta_j \eta\\
&& +\biggl[(\rho_j'-c_j)-\frac {2(p-3)}{\int R_j^2}\int R_j^p\eta\biggr]
\int R_{jx} \Psi_j \eta=E_1+E_2+E_3+E_4+E_5,
\end{eqnarray*}
where we have used
\begin{equation*}
-\int (\eta_{xx}+pR_j^{p-1} \eta)_x \Theta_j \eta-\frac {c_j}2 \int \Theta_{jx} \eta^2
=-H_j(\eta,\eta)+\frac 12 \int \Phi_{jxx} \eta^2-\frac{p}2 \Theta_j(0)\int (R_j^{p-1})_x \eta^2,
\end{equation*}
 since $\frac {\partial}{\partial x} \Psi_j=\frac {\partial}{\partial x} \Theta_j=\Phi_j$
 and  $\Theta_j(t)=\Theta_j(0)+\Psi_j(t)$.

From this identity, we claim
\begin{claim}\label{VIR} There exist $A\geq  5$, $\kappa_0>0$ such that, for $\alpha$ small enough, 
and for all $t\geq 0$, 
\begin{equation}
\frac {d}{dt} \mathcal{K}_j(t)\le -\kappa_0 \sqrt{c_j}\int (\eta_x^2 + c_j\eta^2)e^{-\frac{\sqrt{c_j}}{A}|x-\rho_j|}+\frac 1{\kappa_0}e^{-\frac 18\sqrt{c_j} (t+T_c)}.
\end{equation}
\end{claim}

By \eqref{centieme} and $A\geq 5$, we have $\sqrt{c_1(t)}/A<1/4$ and $\sqrt{c_2(t)}/A<\sqrt{c}/4$
and thus we obtain the conclusion of Claim \ref{VIRIEL} from Claim \ref{VIR}.

\medbreak

\noindent\emph{Proof of Claim \ref{VIR}.}
By \cite{MMarchives} Proposition 6 and localization arguments as \cite{MMjams} proof of Proposition 6  (see also proof of Claim \ref{quadra} in the present paper),   there exists $A_0>0$, $\lambda_0>0$, such that if $A>A_0$ then
\begin{equation}\label{positivite}
H_j^*(\eta,\eta)\ge \lambda_0 \, \sqrt{c_j(t)}\int(\eta_x^2 + c_j\eta^2)e^{-\frac{\sqrt{c_j}}{A}|x-\rho_j|}.
\end{equation}
The claim means that all other terms in the previous identity can be absorbed by this term for some $A\geq  5$ for $\alpha$, $c$ small enough up to and error term of size $e^{-\frac 18 \sqrt{c_j}(t+T_c)}$.

First, since $\Phi_{jxx}=\frac {c_j}{A^2} \Phi_j,$ we have
\begin{equation}
\left|\int \Phi_{jxx} \eta^2\right|\le \frac {K}{A^2}c_j^{\frac 32} \int \eta^2 e^{-\frac{\sqrt{c_j}}{A}|x-\rho_j|}\le \frac{\lambda_0}{100} c_j^{\frac 32}\int \eta^2 e^{-\frac{\sqrt{c_j}}{A}|x-\rho_j|},
\end{equation}
for $A$ large enough. Now, $A$ is fixed to such value.

Next, we have from  $\|\eta\|_{H^1}\le K \alpha c^q$ and (\ref{geom1})
\begin{equation*}
|\rho_j'-c_j|\le K\alpha c_j^{\frac 12 -q} \sqrt{g_j(t)} +Kc_j^{-2q} \int(\eta_x^2 + c_j\eta^2)e^{-\frac{\sqrt{c_j}}{A}|x-\rho_j|}
+K e^{-\frac 18\sqrt{c_j}(t+T_c)},
\end{equation*}
$$
\int \Phi_j \eta^2 \leq K\sqrt{c_j} \int \eta^2e^{-\frac{\sqrt{c_j}}{A}|x-\rho_j|}
\leq K \alpha c^{q} \sqrt{g_j}  \quad \text{and}\quad
\int \Phi_j \eta^2 \leq K \alpha^2 c_j^{2q+\frac 12}.
$$
Therefore,
\begin{equation*}
|(\rho_j'-c_j) \int \Phi_j \eta^2|\le K\alpha \sqrt{c_j}\int(\eta_x^2 + c_j\eta^2)e^{-\frac{\sqrt{c_j}}{A}|x-\rho_j|} +K e^{-\frac 18\sqrt{c_j}(t+T_c)}.
\end{equation*}
Thus, for $\alpha$ small enough,  $|(\rho_j'-c_j) \int \Phi_j \eta^2|\le\frac{\lambda_0}{100}\sqrt{c_j}\int(\eta_x^2 + c_j\eta^2)e^{-\frac{\sqrt{c_j}}{A}|x-\rho_j|}+K e^{-\frac 18\sqrt{c_j}(t+T_c)}$.

By $\|\eta\|_{H^1}\le K \alpha c^q$, $|(x-\rho_j)\Phi_j|\le K$, (\ref{geom2}), we have
\begin{equation*}
\left| \frac {c_j'}{4c_j}\int (x-\rho_j)\Phi_j \eta^2\right|\le
K\alpha^2 |(c_j^{2q})'|\le \frac{\lambda_0}{100} \sqrt{c_j}\int(\eta_x^2 + c_j\eta^2)e^{-\frac{\sqrt{c_j}}{A}|x-\rho_j|}+K e^{-\frac 18\sqrt{c_j}t},
\end{equation*}
for $\alpha$ small enough. Thus
\begin{equation*}
E_1\le -\frac{3\lambda_0} 4  \sqrt{c_j(t)}\int(\eta_x^2 + c_j\eta^2)e^{-\frac{\sqrt{c_j}}{A}|x-\rho_j|}.
\end{equation*}
Now, we treat $E_2$. First, all terms containing powers of $R_k$ are controlled by $Ke^{-\frac 18\sqrt{c_j}(t+T_c)}$
since $\Theta_j$ is exponentially small around $\rho_k(t)$. Therefore, we need only estimate terms of the form 
\begin{equation}\label{ouf} r\in \{2,\ldots,p-1\},\quad
\int R_j^{p-r}\eta^{r}(\Theta_j\eta)_x\quad
\mathrm{and}\quad 
\int \eta^p(\Theta_j\eta)_x=\frac p{p+1}\int \eta^{p+1} \Phi_j.
\end{equation}
For the last term in (\ref{ouf}), we  have, arguing as in the proof of (\ref{encoreclaim}),
\begin{eqnarray*}
\left|\int \eta^{p+1} \Phi_j\right|&\le&
K \sqrt{c_j}\int (\alpha^{p-1} c_j \eta^2+ \alpha^{p-5}c_j^{\frac 32-\frac 4{p-1}} \eta^6)e^{-\frac{\sqrt{c_j}}{A}|x-\rho_j|}
\\
&\le&K\alpha^{p-1} \sqrt{c_j}\int(\eta_x^2 + c_j\eta^2)e^{-{\sqrt{c_j}}|x-\rho_j|}\le \frac{\lambda_0}{100} \sqrt{c_j}\int(\eta_x^2 + c_j\eta^2)e^{-\frac{\sqrt{c_j}}{A}|x-\rho_j|},
\end{eqnarray*}
for $\alpha$ small enough. For the first term in (\ref{ouf}), we  integrate by parts, so that
\begin{equation*}
\int R_j^{p-r}\eta^{r}(\Theta_j\eta)_x
=\frac 1 {r+1} \int (r  R_j^{p-r}\Phi_j- (  R_j^{p-r})_x \Theta_j )\eta^{r+1}.
\end{equation*}
Finally, we treat all these terms similarly as in (\ref{encoreclaim}), so that, for $\alpha$ small enough,
\begin{equation*}
|E_2|\le \frac{\lambda_0}{100} \sqrt{c_j}\int(\eta_x^2 + c_j\eta^2)e^{-\frac{\sqrt{c_j}}{A}|x-\rho_j|}+K e^{-\frac 18\sqrt{c_j}(t+T_c)}.
\end{equation*}

For $E_3$, note that by (\ref{vraicp}),
\begin{equation*}
\begin{split}
\left|E_3\right|&=\left|\frac 12 \Theta_j(0) \left((c_j^{2q})'\int Q^2-\frac p2  \int (R_j^{p-1})_x \eta^2+ (\rho_j'-c_j)\int \eta R_{jx}\right)\right|\\
&\leq K \left|\frac {c_j'}{c_j} \int \eta \tilde R_j\right| +K\int |R_{jx}| (|R_j|^{p-3}\eta^3
+|\eta|^p).
\end{split}\end{equation*}
For the first term, 
 since  $\|\tilde R_j\|_{L^2}\le K c^q$ and
$\|\eta\|_{L^2}\le K \alpha c^q$, we have
$\left|\int \tilde R_j\Theta_j\eta\right|\le K   \alpha c^{2q},
$
and so by (\ref{geom2}),  and  arguing for the other term as for $E_2$,  
for $\alpha$ small enough, we get
\begin{equation*}
|E_3|\le \frac{\lambda_0}{100} \sqrt{c_j}\int(\eta_x^2 + c_j\eta^2)e^{-\frac{\sqrt{c_j}}{A}|x-\rho_j|}+K e^{-\frac 18\sqrt{c_j}(t+T_c)}.
\end{equation*}

For the first term in $E_4$, we proceed as for $E_3$.
The last two terms in $E_4$ are controled by $K e^{-\frac 18\sqrt{c_j}(t+T_c)}$
since they contain products of  exponentially decaying functions centered around $\rho_j(t)$ and $\rho_k(t)$.
Thus,
\begin{equation*}
|E_4|\le \frac{\lambda_0}{100} \sqrt{c_j}\int(\eta_x^2 + c_j\eta^2)e^{-\frac{\sqrt{c_j}}{A}|x-\rho_j|}+K e^{-\frac 18\sqrt{c_j}(t+T_c)}.
\end{equation*}

Since $|\Psi_j|\le K $, $\|R_{jx}\|_{L^2}\le K c_j^{q+\frac 12}$ and
$\|\eta\|_{L^2}\le K \alpha c^{q}$, we have
$\left|\int R_{jx} \Psi_j \eta\right|\le K  \alpha c^{2q+\frac 12},
$
and then we obtain by (\ref{geom4})
\begin{equation*}
|E_5|\le \frac{\lambda_0}{100} \sqrt{c_j}\int(\eta_x^2 + c_j\eta^2)e^{-\frac{\sqrt{c_j}}{A}|x-\rho_j|}+K e^{-\frac 18\sqrt{c_j}(t+T_c)},
\end{equation*}
for $\alpha$ small enough.
Thus, the proof of Claim \ref{VIR} is complete.
\hfill$\Box$

\subsection{Proof of  Claim  \ref{MONOETA} -- Monotonicity results on $\eta(t)$}
In this appendix, we prove   monotonicity results for quantities defined in $\eta(t)$.
Claim \ref{MONOETA} is a direct consequence of Claim \ref{NEW} below for $x_0=0$.
Recall that $\psi(x)$ is defined by  \eqref{surphiff}.

For  $0\leq t_0\leq t$, $x_0\geq 0$, $j=1,2$, let
\begin{equation*}
\begin{split}
& \mathcal{M}_j(t)=\int \eta^2 \psi_j ,\\
&\mathcal{E}_j(t)=\int  \left[\frac 12 \eta_x^2
-\frac 1{p{+}1}\left( (R_1{+}R_2{+}\eta)^{p+1} {-}(p{+}1) R_1^p \eta {-} (p{+}1) R_2^p \eta {-}(R_1{+}R_2)^{p+1}\right)\right] \psi_j,
\end{split}
\end{equation*}
where $\psi_1(x)=\psi(\tilde x)$,
$\tilde x= x-\rho_1(t)+x_0+\frac 12 (t-t_0)$,
and $\psi_2(x)=\psi(\sqrt{c}\tilde x_c)$, $\tilde x_c= x-\rho_2(t)+x_0+ \frac c2 (t-t_0) $.
 
For $0\leq t\leq t_0$, $x_0\geq 0$, $\bar x= x-\rho_1(t)-x_0-\frac 12 (t_0-t)$, let
\begin{equation*}
\begin{split}
&  {\mathcal{M}}_0(t)=\int \eta^2(t) \psi( \bar x)  ,\\
&  {\mathcal{E}}_0(t)=\int  \left[\frac 12 \eta_x^2
-\frac 1{p{+}1}\left( (R_1{+}R_2{+}\eta)^{p+1} {-}(p{+}1) R_1^p \eta {-} (p{+}1) R_2^p \eta {-}(R_1{+}R_2)^{p+1}\right)\right]
\psi(\bar x).
\end{split}
\end{equation*}
For $t_0\geq \frac 8 c T_c$, let $\bar t_0=\frac c8 t_0$, and for $\bar t_0\leq t\leq t_0$,
let $\bar x_c= x-\rho_2(t)-\sigma_0 (t_0-t)$,
where $\sigma_0= {(\rho_1(\bar t_0)-\rho_2(\bar t_0))}/{(t_0-\bar t_0)}$; note that
$\frac c {16} \leq \sigma_0 \leq \frac c4$ by \eqref{quatorze}. Let
\begin{equation*}
\begin{split}
&  {\mathcal{M}}_{1,2}(t)=\int \eta^2(t ) \psi(\sqrt{c}\,\bar x_c)  ,\\
&  {\mathcal{E}}_{1,2}(t)=\int  \left[\frac 12 \eta_x^2
-\frac 1{p{+}1}\left( (R_1{+}R_2{+}\eta)^{p+1} {-}(p{+}1) R_1^p \eta {-} (p{+}1) R_2^p \eta {-}(R_1{+}R_2)^{p+1}\right)\right]
\psi(\sqrt{c}\, \bar x_c).
\end{split}
\end{equation*} 
 
\begin{claim}\label{NEW}
Let $x_0>0$, $t_0>0$.
(i)
For all $t\geq t_0 $,
\begin{equation*}
\begin{split}
	& \frac d{dt}\left( c_1^{2q}(t) \int Q^2+ \mathcal{M}_1(t) \right) 
	   \leq K e^{-\frac 1{16}(t-t_0+x_0)} g_1(t) + Ke^{-\frac 1 {32} \sqrt{c} (t+T_c)},\\
	& \frac d{dt}\left(-\frac {2q}{2q+1} c_1^{2q+1}(t) \int Q^2+ 2 \mathcal{E}_1(t) 
	+\frac 1{100} \left( c_1^{2q}(t) \int Q^2+ \mathcal{M}_1(t) \right)
	\right) 
	   \\ &\qquad \leq K e^{-\frac 1{16}(t-t_0+x_0)} g_1(t)+ Ke^{-\frac 1 {32} \sqrt{c} (t+T_c)}.\\
	   & \frac d{dt}\left( \left(c_1^{2q}(t)+c_2^{2q}(t)\right) \int Q^2+ \mathcal{M}_2(t) \right) 
	   \leq K e^{-\frac {c\sqrt{c}}{16}    (t-t_0)} e^{-\frac {\sqrt{c}}{16}  x_0} \sqrt{c}\, g_2(t) + Ke^{-\frac 1 {32} \sqrt{c} (t+T_c)},\\
	& \frac d{dt}\left(-\frac {2q}{2q{+}1} \left(c_1^{2q+1}(t){+}c_2^{2q+1}(t)\right) \int Q^2+ 2 \mathcal{E}_2(t) 
	+\frac c{100} \left( \left(c_1^{2q}(t){+}c_2^{2q}(t)\right) \int Q^2+ \mathcal{M}_2(t) \right)
	\right) \\ &\qquad
	   \leq K e^{-\frac {c\sqrt{c}}{16}   (t-t_0)  } e^{-\frac {\sqrt{c}}{16}  x_0} c^{\frac 32} g_2(t)+ Ke^{-\frac 1 {32} \sqrt{c} (t+T_c)}.
\end{split}
\end{equation*}
(ii) For all  $0\leq t\leq t_0$, 
\begin{equation*}
\begin{split}
  \frac d{dt} \mathcal{M}_0(t)
& \leq K e^{-\frac 1{16}( t_0-t +x_0)} g_1(t)
  + K e^{-\frac 1{16}( t_0-t +x_0)} e^{-\frac{\sqrt{c}}{32} (t+T_c)},
\end{split}
\end{equation*}
\begin{equation*}
\begin{split}
  \frac d{dt} \left(\mathcal{E}_0(t)+\frac 1{100} \mathcal{M}_0(t)\right)
& \leq K e^{-\frac 1{16}( t_0-t +x_0)} g_1(t)
 + K e^{-\frac 1{16}( t_0-t +x_0)} e^{-\frac{\sqrt{c}}{32} (t+T_c)}.
\end{split}
\end{equation*}
(iii) For all $\bar t_0\leq t\leq t_0$, 
\begin{equation*}
\begin{split}
	   & \frac d{dt}\left(  c_1^{2q}(t)  \int Q^2+  {\mathcal{M}}_{1,2}(t) \right) 
	   \leq K e^{-\frac {c\sqrt{c}}{16}    (t-t_0)} \sqrt{c}\, g_2(t) 
	  + K e^{-\frac {\sqrt{c}}{16} (t-\bar t_0)} g_1(t)
	   + Ke^{-\frac 1 {32} \sqrt{c} (t+T_c)},\\
	& \frac d{dt}\left(-\frac {2q}{2q+1}   c_1^{2q+1}(t)  \int Q^2+ 2  {\mathcal{E}}_{1,2}(t)
	+\frac c{100} \left(  c_1^{2q}(t)  \int Q^2+  {\mathcal{M}}_{1,2}(t) \right)
	 \right)
	\\ &\hskip 3cm
	\leq K e^{-\frac {c\sqrt{c}}{16}   (t-t_0)  } c^{\frac 32} g_2(t)
	+ K e^{-\frac {\sqrt{c}}{16} (t-\bar t_0)} g_1(t)
	+ Ke^{-\frac 1 {32} \sqrt{c} (t+T_c)}.
\end{split}
\end{equation*}
\end{claim}
Note that the monotonicity results on $\mathcal{E}_{j}(t)$ requires the addition 
of some quantity related to $\mathcal{M}_j(t)$ (here the constant $\frac 1{100}$ is somewhat
arbitrary, any small fixed positive constant would work). See also Lemma 1 in \cite{Martel}
for similar calculations in $u(t)$.

\medskip

\noindent\emph{Proof of Claim \ref{NEW}.} We prove the part of Claim \ref{NEW}Ê
concerning $\mathcal{M}_1(t)$, $\mathcal{E}_1(t)$, $\mathcal{M}_2(t)$, $\mathcal{E}_2(t)$.
The rest is proved similarly.

Let $\tilde x=x-\rho_1(t)+\sigma (t-t_0) +x_0$, where $0<\sigma\leq \frac 12$, $x_0>0$, $t\geq t_0\geq 0$.
First, we compute $\frac {d}{dt} \int \eta^2(t) \psi(\tilde x)$ 
using the equation of $\eta(t)$ (\eqref{eqeta}).
\begin{equation*}
\begin{split}
\frac 12 \frac d {dt} \int \eta^2 \psi(\tilde x) & = \int \eta_t \eta \psi(\tilde x)
 + (\sigma -\rho_1'(t)) \int \eta^2 \psi'(\tilde x)\\
 & = - \int \eta_{xxx} \eta \psi(\tilde x) - \int ((R_1+R_2+\eta)^p - R_1^p-R_2^p)_x \eta \psi(\tilde x)
 -\frac {c_1'}{2c_1} \int \tilde R_1 \eta \psi(\tilde x) 
 \\ &\quad - \frac {c_2'}{2c_2} \int \tilde R_2 \eta \psi(\tilde x)
 + (\rho_1'-c_1) \int R_{1x} \eta \psi(\tilde x) + (\rho_2'-c_2) \int R_{2x} \eta \psi(\tilde x)
\\ &\quad  + (\sigma -\rho_1'(t)) \int \eta^2 \psi'(\tilde x) .
\end{split}
\end{equation*}
First, $-\int \eta_{xxx} \eta \psi'(\tilde x) =-\frac 32 \int \eta_x^2 \psi'(\tilde x)
+\frac 12 \int \eta^2 \psi'''(\tilde x)$ and so by
$\psi'''\leq \frac 14 \psi'$, $\sigma \leq \frac 12$ and $\rho_1'(t)\geq \frac 34$,
we obtain similarly as for $u(t)$ in the proof of Claim \ref{MONOTONIE},
\begin{equation*}
-\int \eta_{xxx} \eta \psi(\tilde x)  
 + (\sigma -\rho_1'(t)) \int \eta^2 \psi'(\tilde x) 
 \leq -\frac 32 \int \eta_x^2 \psi'(\tilde x)
-\frac 18\int \eta^2 \psi' (\tilde x).
\end{equation*}
Note that by the decay properties of $R_2$ and $\psi$, we have (see e.g. Claim \ref{prems})
\begin{equation*}
\left | \int \tilde R_2 \eta \psi(\tilde x)\right|
+ \left | \int R_{2x} \eta \psi(\tilde x)\right|
+\left| \int ((R_1+R_2+\eta)^p - (R_1+\eta)^p)_x \eta \psi(\tilde x)\right|
\leq K e^{-\frac 1{32} \sqrt{c} (t+T_c)}.
\end{equation*}
Thus,
\begin{equation*}
\begin{split} 
\frac 12 \frac d {dt} \int \eta^2 \psi(\tilde x) & \leq 
-\frac 32 \int \eta_x^2 \psi'(\tilde x)
-\frac 18\int \eta^2 \psi' (\tilde x)  
  - \int ((R_1 +\eta)^p - R_1^p )_x \eta \psi(\tilde x)\\ &
 -\frac {c_1'}{2 c_1} \int \tilde R_1 \eta \psi(\tilde x) 
  + (\rho_1'-c_1) \int R_{1x} \eta \psi(\tilde x) +K e^{-\frac 1{32} \sqrt{c} (t+T_c)}.
\end{split}
\end{equation*}
Note that 
\begin{equation*}
\begin{split}
&\int \left[(R_1+\eta)^p-R_1^p-pR_1^{p-1}\eta \right] R_{1x} 
=\int  (R_1+\eta)^p ((R_1+\eta)_x -\eta_x)-\int (R_1^p)_x\eta
\\ &=-\int (R_1+\eta)^p \eta_x -(R_1^p)_x \eta = \int ((R_1+\eta)^p -R_1^p)_x\eta,
\end{split}
\end{equation*}
and from  \eqref{vraicp},
\begin{equation*}\begin{split}
&  \frac 12(c_1^{2q})'\int Q^2\leq 
  -(\rho_1'-c_1) \int \eta R_{1x} +\frac {c_1'}{2c_1} \int \eta \tilde R_1  \\&
   +\int  ((R_1+\eta)^p -R_1^p)_x\eta
+ K e^{-\frac 1{8}\sqrt{c}(t+T_c)}.\end{split}\end{equation*} 
We obtain
\begin{equation*}
\begin{split}
&\frac 12 \frac {d}{dt}  \left(c_1^{2q}(t)\int Q^2+ \mathcal{M}_1(t)\right)
 \leq  -\frac 32 \int \eta_x^2 \psi'(\tilde x)
-\frac 18\int \eta^2 \psi' (\tilde x) + K e^{-\frac 1{8}\sqrt{c}(t+T_c)}\\ &
+ \int ((R_1 +\eta)^p - R_1^p )_x \eta (1-\psi(\tilde x))
 +\frac {c_1'}{2c_1} \int \tilde R_1 \eta (1-\psi(\tilde x) )
    - (\rho_1'-c_1) \int R_{1x} \eta (1-\psi(\tilde x)) \\
  &+ K e^{-\frac 18 (\sigma (t-t_0) +x_0)} g_1(t)
  +K e^{-\frac 1{32}\sqrt{c}(t+T_c)}.
  \end{split}
\end{equation*}
Therefore,
 by  $(1-\psi(\tilde x)) e^{-\sqrt{c_1(t)}(x-\rho_1(t))}\leq K e^{-\frac 18 (\sigma (t-t_0) +x_0)}$,
\begin{equation}\label{NUMM}
\begin{split}
\frac 12 \frac {d}{dt}  \left(c_1^{2q}(t)\int Q^2+ \mathcal{M}_1(t)\right)
& \leq  -\frac 32 \int \eta_x^2 \psi'(\tilde x)
-\frac 18\int \eta^2 \psi' (\tilde x) \\ &
 +K e^{-\frac 18 (\sigma (t-t_0) +x_0)} g_1(t)
  +K e^{-\frac 1{32}\sqrt{c}(t+T_c)}.
  \end{split}
\end{equation}

The proof is similar for $\mathcal{E}_1(t)$, up to some additional algebraic cancellations. First,
\begin{equation*}
\begin{split}
\frac d {d t}\mathcal{E}_1&=
\int \eta_{tx} \eta_x \psi(\tilde x) - \int \eta_t ((R_1+R_2+\eta)^p-R_1^p-R_2^p)  \psi(\tilde x)\\
& - \int ((R_1+R_2+\eta)^p -p R_1^{p-1}\eta - (R_1+R_2)^p) R_{1t} \psi(\tilde x)\\
& - \int ((R_1+R_2+\eta)^p -p R_2^{p-1} \eta- (R_1+R_2)^p) R_{2t} \psi(\tilde x)\\
& + (\sigma -\rho_1'(t)) \int \Bigl[\tfrac 12 \eta_x^2 - \tfrac 1{p+1}  ((R_1+R_2+\eta)^{p+1}
-(p+1) R_1^p \eta \\ & -(p+1) R_2^p \eta - (R_1+R_2)^{p+1})\Bigr] \psi'(\tilde x).
\end{split}
\end{equation*}
All terms containing $R_2$, $R_{2t}$ are controlled by $K e^{-\frac 1{8}\sqrt{c}(t+T_c)}$.
Note that $R_{1t}=\frac {c_1'}{2c_1} \tilde R_1 - \rho_1' R_{1x}$. Thus,
\begin{equation*}
\begin{split}
\frac d {d t}\mathcal{E}_1&\leq 
\int \eta_{t} (-\eta_{xx} - ((R_1+R_2+\eta)^p-R_1^p-R_2^p)  \psi(\tilde x)
-\int \eta_t \eta_x \psi'(\tilde x) \\
& - \int ((R_1 +\eta)^p -p R_1^{p-1}\eta - R_1^p) (\frac {c_1'}{2c_1} \tilde R_1 - \rho_1' R_{1x}) \psi(\tilde x)\\
& + (\sigma -\rho_1'(t)) \int \Bigl[\tfrac 12 \eta_x^2 - \tfrac 1{p+1}  ((R_1 +\eta)^{p+1}
-(p+1) R_1^p \eta    - R_1^{p+1} )\Bigr]\psi'(\tilde x)+K e^{-\frac 1{8}\sqrt{c}(t+T_c)}.
\end{split}
\end{equation*}
Thus, using the equation of $\eta$, we get
\begin{equation*}
\begin{split}
\frac d {d t}\mathcal{E}_1 &\leq
-\frac 12 \int  (-\eta_{xx} - ((R_1+R_2+\eta)^p-R_1^p-R_2^p)^2 \psi'(\tilde x) \\
& -\frac {c_1'}{2c_1}\int  \tilde R_1 (-\eta_{xx} - ((R_1 +\eta)^p - R_1^p))\psi(\tilde x)\\
& + (\rho_1'-c_1) \int R_{1x} (-\eta_{xx} - ((R_1 +\eta)^p - R_1^p))\psi(\tilde x)\\
&+\int (\eta_{xx}+(R_1+\eta)^p -R_1^p)_x \eta_x \psi'(\tilde x) +\frac {c_1'}{2c_1} \int \tilde R_1 \eta_x \psi'(\tilde x)
-(\rho_1' -c_1)\int R_{1x} \eta_x \psi'(\tilde x)\\
& -\int ((R_1+\eta)^p - p R_1^{p-1} \eta -R_1^p) (\frac {c_1'}{2c_1} \tilde R_1 - \rho_1' R_{1x})\psi(\tilde x)\\
& + (\sigma -\rho_1'(t)) \int \Bigl(\tfrac 12 \eta_x^2 - \tfrac 1{p+1}  ((R_1 +\eta)^{p+1}
-(p+1) R_1^p \eta    - R_1^{p+1} )\Bigr) \psi'(\tilde x)
+K e^{-\frac 1{8}\sqrt{c}(t+T_c)}.
\end{split}
\end{equation*}
First, terms containing $\psi'(\tilde x)$ and $R_1$ are all controlled by 
$K e^{-\frac 18 (\sigma (t-t_0) +x_0)} g_1(t)$.
Second, we note that
$$\int \eta_{xxx}\eta_x \psi'(\tilde x)+p \int \eta_x^2\eta^{p-1} \psi'(\tilde x) +(\sigma - \rho_1')\int \tfrac 12 \eta_x^2 \psi'(\tilde x)
\leq \frac 12 \int \eta_x^2 \psi'''(\tilde x)-\frac 1{16} \int \eta_x^2 \psi'(\tilde x)\leq 0,$$ 
by
$\psi'''\leq \frac 1{16} \psi'$, $\frac 14 \leq \rho_1'-\sigma \leq 2$
and the fact that $p\int |\eta_x^2\eta^{p-1}| \psi'(\tilde x)
\leq K \|\eta\|_{H^1}^{p-1} \int \eta_x^2\psi'(\tilde x)\leq \frac 1{16} \int \eta_x^2 \psi'(\tilde x)$
for $\alpha$ small enough.

Therefore, 
\begin{equation*}
\begin{split}
\frac d {d t}\mathcal{E}_1 &\leq
\frac {2}{p+1} \int |\eta|^{p+1} \psi'(\tilde x)
-\frac {c_1'}{2c_1}\int  \tilde R_1 (-\eta_{xx} + c_1\eta -p R_1^{p-1}\eta)\psi(\tilde x)\\
& + (\rho_1'-c_1) \int R_{1x} (-\eta_{xx}  + c_1\eta -p R_1^{p-1}\eta)\psi(\tilde x)
+ c_1 \int R_{1x} ((R_1+\eta)^p - p R_1^{p-1} \eta -R_1^p)\psi\\
&+\frac {c_1'}{2c_1} c_1 \int \tilde R_1 \eta\psi - c_1 (\rho_1'-c_1) \int \eta R_{1x} \psi
+ K e^{-\frac 18 (\sigma (t-t_0) +x_0)} g_1(t)
+K e^{-\frac 1{8}\sqrt{c}(t+T_c)}.
\end{split}
\end{equation*}
Using the fact that $\mathcal{L}\tilde Q= -2 Q$ and $\mathcal{L}Q'=0$, $(\tilde R_1 \psi)_{xx}
=\tilde R_{1xx} \psi+ 2 \tilde R_{1x} \psi' + \tilde R_1 \psi''$,
we obtain:
\begin{equation*}
\begin{split}
\frac d {d t}\mathcal{E}_1 &\leq
\frac {2}{p+1} \int |\eta|^{p+1} \psi'(\tilde x)+
\frac {c_1'}{2c_1} c_1 \int \tilde R_1 \eta\psi - c_1 (\rho_1'-c_1) \int \eta R_{1x} \psi
\\ &
+ c_1 \int R_{1x} ((R_1+\eta)^p - p R_1^{p-1} \eta -R_1^p)\psi
+K e^{-\frac 18 (\sigma (t-t_0) +x_0)} g_1(t)
+K e^{-\frac 1{8}\sqrt{c}(t+T_c)}.
\end{split}
\end{equation*}
Using \eqref{vraicp} and $(1-\psi(\tilde x)) e^{-\sqrt{c_1(t)}(x-\rho_1(t))}\leq K e^{-\frac 18 (\sigma (t-t_0) +x_0)}$, we obtain
$$
\frac d {d t}\mathcal{E}_1\leq
\frac 12 c_1 (c_1^{2q})' \int Q^2 +\frac {2}{p+1} \int |\eta|^{p+1} \psi'(\tilde x)+
K e^{-\frac 18 (\sigma (t-t_0) +x_0)} g_1(t)+K e^{-\frac 1{8}\sqrt{c}(t+T_c)}.
$$
Since $ c_1 (c_1^{2q})'=\frac {2q}{2q+1} (c_1^{2q+1})'$, 
and $\frac {2}{p+1} \int |\eta|^{p+1} \psi'(\tilde x)\leq K |\eta\|_{H^1}^{p-1} \int \eta^2\psi'(\tilde x) \leq \frac 1{400} \int \eta^2 \psi'(\tilde x)$ for $\alpha$ small enough,
by \eqref{NUMM}, we obtain:
\begin{equation*}
\begin{split}&
\frac d{dt}\left(-\frac {2q}{2q+1} c_1^{2q+1}(t) \int Q^2+ 2 \mathcal{E}_1(t) 
	+\frac 1{100} \left( c_1^{2q}(t) \int Q^2+ \mathcal{M}_1(t) \right)
	\right) 
	   \\ &\qquad \leq K e^{-\frac 1{8}(\sigma(t-t_0)+x_0)} g_1(t)+ Ke^{-\frac 1 {32} \sqrt{c} (t+T_c)}.
	   \end{split}
\end{equation*}
Thus the claim is proved for 
 $\mathcal{E}_1(t)$. 

For $\mathcal{M}_2(t)$ and $\mathcal{E}_2(t)$, the proof is the same.
For example, we have, for $t\geq t_0$,
\begin{equation*}
\begin{split}
\frac 12 \frac {d}{dt} 
\mathcal{M}_2(t)
 &\leq  \int ((R_1 +\eta)^p - R_1^p )_x \eta (1-\psi(\sqrt{c} \tilde x_c))
+\int ((R_2 +\eta)^p - R_2^p )_x \eta (1-\psi(\sqrt{c} \tilde x_c))
 \\& +\frac {c_1'}{2c_1} \int \tilde R_1 \eta (1-\psi(\sqrt{c} \tilde x_c) )
 +\frac {c_2'}{2c_2} \int \tilde R_2 \eta (1-\psi(\sqrt{c} \tilde x_c) )
  +K e^{-\frac 1{8}\sqrt{c}(t+T_c)}\\ & - (\rho_1'-c_1) \int R_{1x} \eta (1-\psi(\sqrt{c} \tilde x_c))
  - (\rho_2'-c_2) \int R_{2x} \eta (1-\psi(\sqrt{c} \tilde x_c))
  \\
  &\leq K e^{-\frac {\sqrt{c}} 8 (\sigma  {c}\,(t-t_0) +x_0)} \sqrt{c} \, g_2(t)
  +K e^{-\frac 1{32}\sqrt{c}(t+T_c)},
\end{split}
\end{equation*}
by Claim \ref{params} (as in the proof of Claim \ref{VIR}),
  $  R_1 (1-\psi(\sqrt{c}\tilde x_c)) \leq K e^{-\frac 1{32}\sqrt{c}(t+T_c)}$
and $ R_2 (1-\psi(\sqrt{c}\tilde x_c))\leq K c^{\frac 1{p-1}}  e^{-\frac {\sqrt{c}} 8 (\sigma  {c}\,(t-t_0) +x_0)}$.
\quad
$\Box$

\subsection{Proof of Lemma \ref{OdeT}}
We follow similar steps as in the proof of \eqref{pourstab}, using monotonicity arguments
on $\eta(t)$ instead of $u(t)$.

\medskip

\emph{(i) Estimate in the region $x>\rho_1(t)+x_0$.}
We claim, for all $x_0>0$:
\begin{equation}\label{petiti}\begin{split}
&\int_0^{+\infty} \int (\eta_x^2(t,x)+\eta^2(t,x)) \psi(x-\rho_1(t)-x_0) dx dt
\\ & \leq  K( \alpha^2 c^{2 q+1} + e^{-\frac 74 {x_0} } + \exp(-2 c^{-r})) 
+ K \int_{x>x_0} x^2u^2(0,x) dx.
\end{split}\end{equation}

Let us prove \eqref{petiti}.
By Claim \ref{NEW}, using the estimate on $\frac d{dt}\mathcal{M}_0(t)$
and integrating between $0$ and $t_0$, we get the estimate:
\begin{equation}\label{abve}\begin{split}
\int \eta^2(t_0,x) \psi(x-\rho_1(t_0)-x_0) dx  & \leq \int \eta^2(0,x) \psi(x-\rho_1(0)-x_0-\tfrac 12 t_0) dx\\
&+\int_0^{t_0} e^{-\frac 1{16} ( t_0-t +x_0)} g_1(t) dt 
+ K e^{ -\frac{\sqrt{c}}{32} (t_0+T_c)-\frac 1{16} x_0}.
\end{split}\end{equation}
Note that by Fubini Theorem,   $|\rho_1(0)|\leq 1$, the expression of $\psi$ \eqref{surphiff},
and $\|\eta(0,x)\|^2\leq K \alpha^2 c^{2q+1}$,
\begin{equation*}\begin{split}
 \int_0^{+\infty} \int \eta^2(0,x) \psi(x-\rho_1(0)-x_0-\tfrac 12 t) dxdt &
=2 \int \eta^2(0,x+\rho_1(0)+x_0) \left(\int_{-\infty}^{x} \psi\right) dx\\
&\leq C \alpha^2 c^{2q+1} +  K \int_{x>x_0} x \eta^2(0,x) dx.
\end{split}\end{equation*}
Next, by Lemma \ref{INTETA} and $x_0\geq 0$,
\begin{equation*}\begin{split}
\int_0^{+\infty} \int_0^{t_0} e^{-\frac 1 {16} ((t-t_0) +x_0)} g_1(t) dt dt_0
& \leq \int_0^{+\infty} g_1(t) \Biggl(\int_t^{+\infty} e^{-\frac 1 {16} (t-t_0)  } dt_0\Biggr)dt\\&
\leq  16  \int_0^{+\infty} g_1(t) dt 
\leq K (\alpha^2 c^{2q+1} + \exp(-2 c^{-r})).
\end{split}\end{equation*}
Finally, by $\int_0^{+\infty}
 e^{ -\frac{\sqrt{c}}{32} (t+T_c)-\frac 1{16} x_0} dt \leq \exp(-2 c^{-r})$ and integrating
\eqref{abve} for $t_0\in [0,+\infty)$, we obtain
\begin{equation}\label{abveint}\begin{split}
& \int_0^{+\infty} 
\int \eta^2(t ,x) \psi(x-\rho_1(t)-x_0) dx dt   \leq 
 \int_{x>x_0} x \eta^2(0,x) dx
+ K (\alpha^2 c^{2q+1} + \exp(-2 c^{-r})).
\end{split}\end{equation}

To control $\int_0^{+\infty} \int  \eta_x^2(t,x) \psi(x-\rho_1(t)-x_0) dx dt$,
we use the same argument on $\mathcal{E}_0(t)$.
First, we claim as a consequence of a standard argument (Kato's identity) applied to the equation of
$\eta(t)$ that there exists $0< \bar t\leq 1$, such that
\begin{equation}\label{gngn}
\int_{x>0} x \, (\eta_x^2+\eta^2) (\bar t,x) dx \leq K \alpha^2 c^{2q+1} + K \int_{x>x_0} x^2\eta^2(0,x) dx.
\end{equation}
 Indeed, 
let $\lambda>0$ to be fixed large enough and
$I(t)=\int \eta^2(t,x) \left(\int_{-\infty}^{x-\lambda t -x_0} (\int_{-\infty}^s \psi) ds\right) dx.$
Then, by the
equation of $\eta$ \eqref{eqeta}, for some $K_0>0$,
\begin{equation*}\begin{split}
\frac 12 I'(t) &\leq
-\frac 32 \int \eta_x^2  \int_{-\infty}^{x-\lambda t -x_0}  \psi 
-\lambda \int \eta^2  \int_{-\infty}^{x-\lambda t -x_0}  \psi\\
& 
+\int \eta^2 \psi'(x-\lambda t -x_0) +  K_0 \int  \eta^2  \int_{-\infty}^{x-\lambda t -x_0}  \psi
 \leq - \frac 12 \int (\eta_x^2+\eta^2)  \int_{-\infty}^{x-\lambda t -x_0}  \psi,
\end{split}\end{equation*}
by $\psi' \leq \frac 14 \int_{-\infty}^x \psi$ (from \eqref{surphiff}) and choosing  $\lambda > K_0+1$.
Thus, $\int_0^1 \int ( \eta_x^2+\eta^2)
\int_{-\infty}^{x-\lambda t -x_0}  \psi dx dt\leq 2 I(0)$,
and there exists $\bar t\in [0,1]$ such that
\begin{equation}\label{estgrad}\begin{split}
 & \int (\eta_x^2+\eta^2)(\bar t)  \int_{-\infty}^{x -\rho_1(\bar t)-x_0}  \psi dx
 \leq K \int (\eta_x^2+\eta^2)(\bar t)  \int_{-\infty}^{x-\lambda \bar t -x_0}  \psi dx
\\ &\leq K I(0)
\leq K \alpha^2 c^{2q+1} + K \int_{x>x_0} x^2\eta^2(0,x) dx.
\end{split}\end{equation}

By   Claim \ref{NEW},  Lemma \ref{INTETA}, \eqref{abveint} and $\|\eta(t)\|_{H^1}\leq K \alpha^2c^{2q+1}$ for $t\in [0,1]$,
we obtain
\begin{equation}\label{abveintX}\begin{split}
& \int_0^{+\infty}
\int \eta_x^2(t ,x) \psi(x-\rho_1(t)-x_0) dx dt  \leq 
\int_0^{\bar t}\int \eta_x^2(t ,x) \psi(x-\rho_1(t)-x_0) dx dt\\ &
 + 2\int_{\bar t}^{+\infty} \mathcal{E}_{0}(t) dt
 + K \int_{\bar t}^{+\infty}\int \eta^2(t ,x) \psi(x-\rho_1(t)-x_0) dx dt\\ &
\leq 
 K   \int (\eta^2_x+\eta^2) (\bar t,x+\rho_1(\bar t)+x_0)  \left(\int_{-\infty}^{x}  \psi\right) dx 
+ K( \alpha^2 c^{2 q+1} +  \exp(-2 c^{-r})).
\end{split}\end{equation}
From \eqref{estgrad}, we obtain, for all $x_0>0$,
\begin{equation}\label{DROITE}\begin{split}
&\int_0^{+\infty} \int (\eta_x^2(t,x)+\eta^2(t,x)) \psi(x-\rho_1(t)-x_0) dx dt
\\ & \leq  K( \alpha^2 c^{2 q+1} +  \exp(-2 c^{-r})) 
+ K \int_{x>x_0} x^2\eta^2(0,x) dx.
\end{split}\end{equation}

Since $\eta(0,x)=u(0,x)-Q_{c_1(0)}(x-\rho_1(0))-Q_{c_2(0)}(x -\rho_2(0))$,
and $\rho_1(0)-\rho_2(0)\geq T_c/4$, \eqref{centieme},
we have
\begin{equation}\label{D122}\begin{split}
& \int_{x>x_0} x^2\eta^2(0,x) dx
  \\ &
\leq K\int_{x>x_0} x^2u^2(0,x) dx+K\int_{x>x_0} x^2(Q_{c_1(t)}^2(x-\rho_1(0))  
+  Q_{c_2(t)}^2(x-\rho_2(0)) )dx\\
& \leq K \int_{x>x_0} x^2u^2(0,x) dx + K ((x_0^2+1) e^{-2\sqrt{c_1(0)} x_0} +  \exp(-2 c^{-r}))
\\ &
\leq K\int_{x>x_0} x^2u^2(0,x) dx + K ( e^{-\frac 74  {x_0} } +  \exp(-2 c^{-r})).
\end{split}\end{equation}
Thus, \eqref{petiti} is proved.

\medskip

\emph{(ii) Estimate of $\tilde g_1(t)$.}
We claim 
\begin{equation}\label{surguntilde}\begin{split}
  \int_0^{+\infty} \tilde g_1(t)  dt
\leq K  ((\alpha c^{q+\frac 12})^{\frac 74} +  \gamma_0(\alpha,c) + \alpha^{-\frac 14}
\exp(-\tfrac 32 c^{-r})),
\end{split}\end{equation}
where $\gamma_0(\alpha,c)= \int_{x>|\ln (\alpha c^{q+\frac 12})|} x^2 u^2(0,x) dx$.
For $x_0>1$ to be fixed later, we observe that
$$0<\psi(x)-\psi(x-x_0)=\int_{x-x_0}^x \psi'(s) ds \leq  K \int_{x-x_0}^x e^{-\frac {|s|}4} ds
\leq K e^{-\frac {|x|}4} e^{ \frac {x_0}4}.$$
Thus,
\begin{equation*}\begin{split}
  \tilde g_1(t)  
\leq K e^{ \frac {x_0} 4}   g_1(t)   + \int (\eta^2_x+\eta^2)(t,x) \psi(x-\rho_1(t)-x_0)dx.
\end{split}\end{equation*}
Therefore, it follows from \eqref{petiti} and Lemma \ref{INTETA} that
\begin{equation*}\begin{split}
  \int_0^{+\infty} \tilde g_1(t)  dt
\leq K \left( e^{ \frac {x_0} 4}   \alpha^2 c^{2q+1} +     e^{-\frac 74{x_0}}\right)+
K \int_{x>x_0} x^2u^2(0,x) dx  +    K e^{ \frac {x_0} 4} \exp(-2 c^{-r})  .
\end{split}\end{equation*}
We obtain \eqref{surguntilde} by choosing $x_0=|\ln (\alpha c^{q+\frac 12})|$.   
 
 \medskip

\emph{(iii) Estimate of $\tilde g_2(t)$.}
We claim
\begin{equation}\label{surgdeuxtilde}
\int_0^{+\infty} \tilde g_2(t)dt
\leq K (\alpha^{\frac 74} c^{\frac 74 q-\frac 18}  + \alpha^2 c^{\frac 32q -\frac 12} + \tfrac 1 c \gamma_0(\alpha,c)+ \alpha^{-\frac 14} \exp(-\tfrac 32 c^{-r}).
\end{equation}
We estimate separately $\int_{0}^{\frac 8 c T_c} \tilde g_2(t) dt$
and $\int_{\frac 8 c T_c}^{+\infty} \tilde g_2(t) dt$.
For the first integral, we use   $\tilde g_2(t)\leq \|\eta(t)\|_{H^1_c}^2\leq
K \alpha^2 c^{2q+1}+ K  \exp(- 2 c^{-r})$ by \eqref{MMnonlinearity}. Thus, by $T_c \leq c^{-\frac 12(1+q)}$,
$$\int_{0}^{\frac 8 c T_c} \tilde g_2(t) dt \leq  \frac Kc T_c \alpha^2 c^{2q+1}
\leq K \alpha^2 c^{\frac 32 q-\frac 12} + K  \exp(- 2 c^{-r}).$$

Now, we use Claim \ref{NEW}. By integration in time on $[\frac c8 t_0,t_0]$ of the
estimates of $\frac d{dt}\mathcal{M}_{1,2}(t)$ and $\frac d{dt}\mathcal{E}_{1,2}(t)$
and using 
$  {\mathcal{E}}_{1,2}(t) + c_1(\frac c8 t_0){\mathcal{M}}_{1,2}(t)\geq \frac 1K \tilde g_2(t)$
(see Claim \ref{quadra}), 
and then   Claim \ref{SURLESC}, we obtain
\begin{equation*}\begin{split}
& \tilde g_2(t_0) \leq
 K\tilde g_1(\tfrac c8 t_0) + K |c_1(t_0)-c_1(\bar t_0)| 
 + K \int_{\frac c8 t_0}^{t_0} \left(e^{-\frac {c\sqrt{c}}{16}   (t_0-t)  } c^{\frac 32} g_2(t)
	+   e^{-\frac {\sqrt{c}}{16} (t-\bar t_0)} g_1(t)\right) dt 
	\\ &+ Ke^{-\frac 1 {32} \sqrt{c} (\frac c8 t_0 t+T_c)}\\
	&\leq
  K (\tilde g_1(\tfrac c8 t_0) + \tilde g_1(t_0))
 + K \int_{\frac c8 t_0}^{t_0} \left(e^{-\frac {c\sqrt{c}}{16}   (t_0-t)  } c^{\frac 32} g_2(t)
	+   e^{-\frac {\sqrt{c}}{16} (t-\bar t_0)} g_1(t)\right) dt 
	+ Ke^{-\frac 1 {32} \sqrt{c} (\bar t_0+T_c)}.
\end{split}\end{equation*}
Integrating in $t_0$ on $[\frac 8 c T_c, +\infty)$, we get ($\bar t=\frac c8 t$) by 
Fubini Theorem,
\begin{equation*}\begin{split}
 \int_{\frac 8 c T_c}^{+\infty} \tilde g_2(t) dt & \leq
 \frac {8K} c \int_{ T_c}^{+\infty} \tilde g_1(\bar t) d\bar t
 + K \int_{\frac 8 c T_c}^{+\infty} \tilde g_1(t) dt  
+ K\int_{T_c}^{+\infty} g_2(t) c^{\frac 32} \int_{t}^{\frac 8c t} e^{-\frac {c\sqrt{c}}{16}   (t_0-t) }dt_0
dt
\\ &
+ K\int_{T_c}^{+\infty} g_1(t) \int_{t}^{\frac 8c t} e^{-\frac {\sqrt{c}}{16}   (t-\frac c8 t_0) } dt_0 dt
+K\exp(-2 c^{-r})
\\Ê& \leq\frac K c \int_0^{+\infty} \tilde g_1(t) dt + K \int_{0}^{+\infty} g_2(t) dt +K\exp(-2 c^{-r}) \\
& \leq\frac K c ( (\alpha c^{q+\frac 12})^{\frac 74}+\gamma_0(\alpha,c))
+ K \alpha^2 c^{2q -\frac 12} +K \alpha^{-\frac 14} \exp(-\tfrac 32 c^{-r}),
\end{split}\end{equation*}
by \eqref{surguntilde} and Lemma \ref{INTETA}.
Thus, \eqref{surgdeuxtilde} is proved.

\medskip

\emph{(iv) Pointwise estimates.}
Now, we claim pointwise estimates, useful for the proof of Lemma \ref{semiconv},
$$
\lim_{t\to +\infty} t (\tilde g_1(t)+\tilde g_2(t)) =0,\quad 
\lim_{t\to +\infty} t \int (\eta_x^2(t,x)+\eta^2(t,x)) \psi(x-\tfrac c {10} t) dx=0,$$
\begin{equation}\label{uniftg1}
\forall t\geq 0, \quad t\, \tilde g_1(t)\leq K\int_0^{+\infty} \tilde g_1(t) dt +K \exp(-2c^{-r}).
\end{equation}
\begin{equation}\label{uniftg2}
t\tilde g_2(t)\leq K \int_0^{+\infty} \tilde g_2(t) dt + \frac K c \int_0^{+\infty} \tilde g_1(t) dt
+  K \exp(-2 c^{-r}).
\end{equation}
We check the estimate for $\tilde g_1(t)$, the other estimates follow from similar arguments.

 Let $t_1\leq t_0\leq t\leq 2t_1$.
Integrating the conclusion of Claim \ref{INTETA} between $t_0$ and $t$, we get
\begin{equation*}
\begin{split}&
    \mathcal{M}_1(t)-  \mathcal{M}_1(t_0)  \leq  
     ( c_1^{2q}(t_0)- c_1^{2q}(t))\int Q^2
  +K\int_{t_0}^t e^{-\frac 1{16} (t'-t_0)} g_1(t')dt'+Ke^{-\frac 1 {32} \sqrt{c} (t_0+T_c)}
	,\\
&  2\mathcal{E}_1(t)-2 \mathcal{E}_1(t_0) 
+\tfrac 1{100}(\mathcal{M}_1(t)-  \mathcal{M}_1(t_0))\leq  K\int_{t_0}^t e^{-\frac 1{16}(t'-t_0) } g_1(t')dt'
  +  Ke^{-\frac 1 {32} \sqrt{c} (t_0+T_c)}
\\ &+  \left[-\frac {2q}{2q+1} (c_1^{2q+1}(t_0)- c_1^{2q+1}(t))  
+\tfrac 1{100} ( c_1^{2q}(t_0)- c_1^{2q}(t))\right] \int Q^2 
		.
\end{split}
\end{equation*}
Since
$\left| c_1^{2q+1}(t)-c_1^{2q+1}(t_0) - \frac {2q+1}{2q}  c_1(t_0)
( c_1^{2q}(t_0)- c_1^{2q}(t))\right|\leq 
K |c_1(t)-c_1(t_0)|^2,
$
we obtain
\begin{equation*}\begin{split}
\mathcal{E}_1(t) + c_1(t_0) \mathcal{M}_1(t)
&\leq \mathcal{E}_1(t_0) + c_1(t_0) \mathcal{M}_1(t_0)+  
K |c_1(t)-c_1(t_0)|^2 \\ &+K\int_{t_0}^t e^{-\frac 1{16} (t'-t_0)} g_1(t')dt'
	+ Ke^{-\frac 1 {32} \sqrt{c} (t_0+T_c)}.
\end{split}\end{equation*}
Thus,
by Claim \ref{SURLESC}, 
\begin{equation*}\begin{split}
\mathcal{E}_1(t) + c_1(t_0) \mathcal{M}_1(t)
&\leq \mathcal{E}_1(t_0) + c_1(t_0) \mathcal{M}_1(t_0)+ K \, \tilde g_1^2(t)+ K \, \tilde g_1^2(t_0)\\
& +K\int_{t_0}^t e^{-\frac 1{16} (t'-t_0)} g_1(t')dt'
	+ Ke^{-\frac 1 {32} \sqrt{c} (t_0+T_c)}.
\end{split}\end{equation*}
We clearly have $\mathcal{E}_1(t_0) + c_1(t_0) \mathcal{M}_1(t_0)\leq K \, \tilde g_1(t_0)$
and by a variant of Claim  \ref{quadra}, there exists
$\kappa_0>0$ such that
$\mathcal{E}_1(t) + c_1(t_0) \mathcal{M}_1(t)
\geq \kappa_0 \int (\eta_x^2(t,x)+\eta^2(t,x)) \psi(\tilde x) dx \geq \frac 1K  \tilde g_1(t)$.
Thus, we obtain
\begin{equation*}\begin{split}
\tilde g_1(t) \leq K \, \tilde g_1(t_0) +  K\int_{t_0}^t e^{-\frac 1{16} (t'-t_0) } g_1(t')dt'
	+ Ke^{-\frac 1 {32} \sqrt{c} (t_0+T_c)}.
\end{split}\end{equation*}
In particular, taking $t=2t_1$ and integrating in $t_0\in [t_1,2t_1]$, we obtain
\begin{equation*}\begin{split}
t_1 \, \tilde g_1(2t_1) & \leq \int_{t_1}^{2t_1} \tilde g_1(t)dt 
+  K\int_{t_1}^{2t_1}\int_{t_0}^{2t_1} e^{-\frac 1{16} (t'-t_0) } g_1(t')dt' dt_0
	+ K\, t_1 e^{-\frac 1 {32} \sqrt{c} (t_1+T_c)}\\
	& \leq K \int_{t_1}^{2t_1} \tilde g_1(t)dt + K\, t_1 e^{-\frac 1 {32} \sqrt{c} (t_1+T_c)}.
\end{split}\end{equation*}
Since $\int_0^{+\infty} \tilde g_1(t) dt<+\infty$, we obtain 
$\lim_{t\to +\infty} t \tilde g_1(t) = 0$ and the estimate on   $t \tilde g_1(t)$.

\subsection{Proof of Lemma \ref{semiconv}}

We claim the following preliminary result on $J_j(t)$ defined in \eqref{defJ}.

\begin{claim}\label{surJ}(i) Equation of $J_j(t)$:
	\begin{equation}\label{surjj}\begin{split}
	& \left| \frac {5-p}{2(p-1)}(\rho_1'(t)-c_1(t) )\int Q^2  + J_1'(t)\right| 
	\leq K \left(g_1(t) + e^{-\frac 1{16}\sqrt{c} (t+T_c)}\right),\\
	 & \left| \frac {5-p}{2(p-1)}(\rho_2'(t)-c_2(t) )\int Q^2  + J_2'(t)\right|
	 \leq K c^{-\frac 1{p-1}} g_1(t) 
	+ K e^{-\frac 1{16}\sqrt{c} (t+T_c)} \\
	&\qquad + K \left(1+ \alpha c^{- \frac 14 q} + c^{-q-\frac 14}   \gamma_0^{\frac 12}(\alpha,c)
	+ \alpha^{-\frac 18} \exp(-\tfrac 12 c^{-r})\right) c^{-2q}   g_2(t) .
	\end{split}\end{equation}
(ii) Estimates on $J_j(t)$:
\begin{equation*}\begin{split}
&|J_1(0)|\leq K ((\alpha c^{q+\frac 12})^{\frac 78}+\gamma_0^{\frac 12}(\alpha,c)+ \exp(-c^{-r})),
\\ &
|J_2(0)|
\leq K (\alpha^{\frac 78} c^{-\frac 14  q+\frac 12}
+ c^{-q+\frac 14} \gamma_0^{\frac 12}(\alpha,c)+ \exp(-\tfrac 12c^{-r})),\\ &
|J_1(t)|+|J_2(t)|\to 0 \quad \text{as $t\to +\infty$}.
\end{split}\end{equation*}
\end{claim}

Assuming this claim, we finish the proof of Lemma \ref{semiconv}.
Integrating \eqref{surjj} in $t$ on $[0,+\infty)$ since 
$|J_1(t)|+|J_2(t)|\to 0$ as $t\to +\infty$, and using Lemma \ref{INTETA}, we obtain
\begin{equation*}\begin{split} 
 \left|\int_0^{+\infty} (\rho_1'(t)-c_1(t)) dt \right|
& \leq  K |J_1(0)| + K \int_0^{+\infty}  g_1(t) dt
+K \exp(-2 c^{-r})
\\ & \leq K((\alpha c^{q+\frac 12})^{\frac 78}+\gamma_0^{\frac 12}(\alpha,c)+ \exp(-c^{-r})),
\end{split}\end{equation*}
\begin{equation*}\begin{split} 
& \left|\int_0^{+\infty} (\rho_2'(t)-c_2(t)) dt \right|
\leq  K |J_2(0)|+ K   c^{-\frac 1{p-1}} \int_0^{+\infty} g_1(t)  dt
+K \exp(-2 c^{-r})
\\ &  + K   \left(1+ \alpha c^{- \frac 14 q} + c^{-q-\frac 14} \gamma_0^{\frac 12}(\alpha,c)+\alpha^{-\frac 18} \exp(-\tfrac 12c^{-r})
	\right)  c^{-2 q } \int_0^{+\infty}  g_2(t) dt
  \\
& \leq K(\alpha^{\frac 78} c^{-\frac 14  q+\frac 12}
 +\alpha^2 c^{-\frac 12}
+ \alpha^3 c^{-\frac 12-\frac 14 q} + (c^{-q+\frac 14} +\alpha^2 c^{-q -\frac 34})
 \gamma_0^{\frac 12}(\alpha,c) + \exp(-\tfrac 14c^{-r})).
\end{split}\end{equation*}

\medskip

\noindent\emph{Proof of Claim \ref{surJ}}
\emph{(i) Equation of $J_j(t)$.}
We  compute the time derivative of $J_j(t)$:
\begin{equation}
	J_j(t)=c_j^{-2q}(t)\int \eta(t,x) \left(\int_{-\infty}^x \tilde R_j(t,x')dx'\right) dx,
\end{equation}
where $\tilde R_j(t,x)=\tilde Q_{c_j(t)}(x-\rho_j(t))$.
The argument is the same as in the proof of Claim \ref{GEOM}.
Let $(j,k)=(1,2)$ or $(j,k)=(2,1)$ and 
$\tilde{\tilde Q}=\frac 2{p-1} \tilde Q + x \tilde Q'$,
$\tilde{\tilde R}_j=c_j^{\frac 1{p-1}}(t) \tilde{\tilde Q}(\sqrt{c_j(t)} (x-\rho_j(t))).$ Then,
\begin{equation*}
\begin{split}
\frac {d}{dt} J_j & =
-2q \frac {c_j' }{c_j } J_j + c_j^{ -2q}   \Biggl(
\int \eta_t \int_{-\infty}^x \tilde R_j 
- \rho_j' \int \eta  \tilde R_{j} 
+\frac {c_j'}{2c_j} \int \eta \int_{-\infty}^x \tilde {\tilde R}_j\Biggr)\\
& =  -2q  \frac {c_j' }{c_j } J_j + c_j^{ -2q}  \Biggl(
-\int  (-\eta_{xx}+c_j \eta -p R_{j}^{p-1} \eta) \tilde R_j -(\rho_j'-c_j)
\int \eta \tilde R_j\\
&\quad + \int ((R_j+R_k+\eta)^p - R_j^p - R_k^p -p R_{j}^{p-1} \eta) \tilde R_j
-\frac {c_1'}{c_1} \int \tilde R_1 \int_{-\infty}^x \tilde R_j
-\frac {c_2'}{c_2} \int \tilde R_2 \int_{-\infty}^x \tilde R_j
\\ & \quad
-(\rho_1'-c_1) \int R_1 \tilde R_j
-(\rho_2'-c_2) \int R_2 \tilde R_j+   \frac {c_j'}{2c_j} \int \eta \int_{-\infty}^x \tilde {\tilde R}_j\Biggr).
\end{split}
\end{equation*}
Note that $\mathcal{L} \tilde Q= - 2 Q$, and thus $\int  (-\eta_{xx}+c_j \eta -p R_{j}^{p-1} \eta) \tilde R_j=-2 \int \eta R_j=0$, 
$\int R_j\tilde R_j=\frac {5-p}{2(p-1)} c_j^{2q} \int Q^2$, $\int \tilde R_j \int_{-\infty}^x \tilde R_j=-c_j^{\frac 2{p-1}-1} \left(\int \tilde Q\right)^2$,
and finally $\left|\int \tilde R_1\int_{-\infty}^x \tilde R_2\right|\leq K c^{-\frac 12 +\frac 1{p-1}}$.
As in the proof of Claim \ref{GEOM}, we obtain
\begin{equation*}
\begin{split}
& \left| J_j'(t) + \frac {5-p} {2(p-1)} (\rho_j'(t)-c_j(t)) \int Q^2 \right|\\
& \leq K \left( \left|\frac {c_j'}{c_j}\right| c^{-q+\frac 14} h_j + |\rho_j'-c_j| c^{-q-\frac 12} \sqrt{g_j}
+c_j^{-2 q} g_j + |c_1'| c_j^{-\frac 1{p-1}} + \left| \frac {c_j'}{c_j}\right| c_j^{-\frac 12}
+ e^{-\frac 18 \sqrt{c} (t+T_c)}\right),
\end{split}
\end{equation*}
where  
$  h_j(t)=\int |\eta(t,x)| \biggl(\int_{-\infty}^x
	e^{-\frac {\sqrt{c_j(t)}}2|x'-\rho_j(t)|} dx'\biggr)dx.
$

Using Claim \ref{GEOM}, we obtain
\begin{equation}\label{firstjj}
\begin{split}
& \left| J_j'(t) + \frac {5-p} {2(p-1)} (\rho_j'(t)-c_j(t)) \int Q^2 \right|\\
& \leq K \left( c_j^{-2q} g_j(t)(1+c^{-q+\frac 34} h_j(t)) + c_j^{-\frac 1{p-1}} g_1(t) 
+ e^{-\frac 1{16} \sqrt{c} (t+T_c)} \right).
\end{split}
\end{equation}

\medskip

\noindent\emph{(ii) Estimates on $h_j$ and $J_j$.} 
We begin with the estimates at $t=0$.
First, by Cauchy-Schwarz inequality, $|\rho_1(0)|\leq 1$,  and for $x_0>1$, we have
($x_+=\max(x,0)$)
\begin{equation*}\begin{split}
\left|J_1(0)\right|^2 &\leq K h_1^2(0) \leq 
K \left(\int |\eta(0)| \psi\right)^2
\leq K  \left(\int   \eta^2(0,x) {(1+x_+^2)} \psi(x) dx\right)\left(\int \frac {\psi(x)}{ 1+x_+^2 }  dx
\right)\\
&\leq K x_0^2 \|\eta(0)\|_{L^2}^2 + K \int_{x\geq x_0} x^2 \eta^2(0,x) dx.
\end{split}\end{equation*}
Choose  $x_0=|\ln(\alpha c^{q+\frac 12})|$ 
so that by $\|\eta(0)\|_{L^2}\leq K \alpha c^{q+\frac 12},$ \eqref{D122}, we obtain
$\left|J_1(0)\right|^2 \leq 
K ((\alpha c^{q+\frac 12})^{\frac 74}+\gamma_0(\alpha,c)+ \exp(-2 c^{-r})).
$

Next, by considering the three space regions
$x<\rho_2(0)$, $\rho_2(0)<x<\rho_1(0)$, $\rho_1(0)<x$, with $\rho_1(0)-\rho_2(0)\leq 2 T_c$, we have
\begin{equation*}\begin{split}
c^{2q-\frac 1{p-1}} \left|J_2(0)\right| & \leq K h_2(0) \leq 
K (c^{-\frac 14} \|\eta(0)\|_{L^2} + T_c^{\frac 12} \|\eta(0)\|_{L^2}+ h_1(0)).
\end{split}\end{equation*}
By $c^{2q-\frac 1{p-1}} = c^{q-\frac 14}$ and
$\|\eta(0)\|_{L^2}\leq K \alpha c^{q+\frac 12},$ we obtain
(using $T_c^{\frac 12}\leq K c^{-\frac 14(1+q)}$),
\begin{equation}\begin{split}
|J_2(0)| & \leq K c^{-q+\frac 14}  (\alpha c^{-\frac 14(1+q)} c^{q+\frac 12} + \alpha^{\frac 78}
c^{\frac {7 q}{8} + \frac 7{16}} + \gamma_0^{\frac 12}(\alpha,c)+\exp(- c^{-r}) ) \\
& \leq K c^{-q+\frac 14}  (\alpha c^{ \frac 34q +\frac 14}  + \alpha^{\frac 78}
c^{\frac {7 q}{8} + \frac 7{16}} + \gamma_0^{\frac 12}(\alpha,c)+\exp(- c^{-r}) )\\
& \leq K  (\alpha^{\frac 78} c^{-\frac q4+\frac 12} + c^{-q+\frac 14}\gamma_0^{\frac 12}(\alpha,c)
+\exp(- c^{-r}) ).
\end{split}\end{equation}

Now, we prove estimates on $h_j(t)$ to be inserted in \eqref{firstjj}.
By the properties of $\psi$, 
\begin{equation*}\begin{split}
h_1^2(t) &\leq 
K  \int   \eta^2(t,x+\rho_1(t)) {(1+x_+^2)} \psi(x) dx
\leq K \int   \eta^2(t,x+\rho_1(t)) \left(\int_{-\infty}^x \int_{-\infty}^s \psi\right) dx.
\end{split}\end{equation*}
By estimate \eqref{abve} and Lemma \ref{INTETA}, for any $t_0, \ x_0>0$, we have
\begin{equation}\label{DFDF} \begin{split}
& \int \eta^2(t_0,x+\rho_1(t_0)) \psi(x -x_0) dx \\
 & \leq \int \eta^2(0,x+\rho_1(0)+\tfrac 12 t_0) \psi(x-x_0) dx
 +\int_0^{t_0} e^{-\frac 1{16} ( t_0-t +x_0)} g_1(t) dt 
+ K e^{ -\frac{\sqrt{c}}{32} (t_0+T_c)-\frac 1{16} x_0}\\
& \leq \int \eta^2(0,x+\rho_1(0)+\tfrac 12 t_0) \psi(x-x_0) dx
+ K e^{-\frac {x_0}{16}} \Bigl(e^{-\frac {t_0}{32}} \alpha^2 c^{2q}
+  \sup_{t\in [\frac 12 t_0,t_0]} g_1(t)+   e^{ -\frac{\sqrt{c}}{32} (t_0+T_c)}\Bigr).
\end{split}\end{equation}
Note that
$\int_0^{+\infty} \int_y^{+\infty} \psi(x-x_0) dx_0 dy =
\int_0^{+\infty} \int_{-\infty}^{x-y} \psi(s') ds' dy =\int_{-\infty}^x \int_{-\infty}^s \psi(s')ds'ds.$
Thus, by Fubini Theorem,
integrating \eqref{DFDF} in $x_0$ on $[y,+\infty)$, for $y\geq 0$ and then
integrating in $y\in [0,+\infty)$, we obtain
\begin{equation*} \begin{split}
& \int \eta^2(t_0,x+\rho_1(t_0))\left(\int_{-\infty}^x \int_{-\infty}^s \psi \right) dx 
  \leq \int \eta^2(0,x+\rho_1(0)+\tfrac 12 t_0) \left(\int_{-\infty}^x \int_{-\infty}^s \psi \right) dx
\\& + K \Bigl(e^{-\frac {t_0}{32}} \alpha^2 c^{2q+1}
+  \sup_{t\in [\frac 12 t_0,t_0]} g_1(t)+   e^{ -\frac{\sqrt{c}}{32} (t_0+T_c)}\Bigr).
\end{split}\end{equation*}
Therefore,  from the assumption on $u(0)$, it follows that 
$J_1(t)\to 0$ as $t\to +\infty$ and for all $t\geq 0$, 
$h_1(t)\leq K$. Estimate (i) for $J_1'$ then follows from \eqref{firstjj}. 

Now, 
we estimate $h_2(t)$. As before, since $\rho_1(t)-\rho_2(t)\leq K (t+T_c)$,
\begin{equation*}\begin{split} 
  h_2(t) & \leq 
 Kc^{-\frac 14} \|\eta(t)\|_{L^2} +K \int_{\rho_2(t)<x<\rho_1(t)} |\eta(t,x)|dx
+ K h_1(t)\\
& \leq 
K \left( c^{-\frac 14} \|\eta(t)\|_{L^2} + c^{-\frac 12} (t+T_c)^{\frac 12}
\sqrt{\tilde g_2(t)} +   h_1(t)\right).
\end{split}\end{equation*}
By Lemma \ref{OdeT}, we obtain that $h_2(t)\to 0$ as $t\to +\infty$.
Moreover, by the estimate on $h_1(t)$, 
 $c^{-\frac 14} \|\eta(t)\|_{L^2}\leq K \alpha c^{q-\frac 14}$, and \eqref{uniftg2}, we have
\begin{equation*}\begin{split} 
  h_2(t) &
  \leq  K\left(1+\alpha c^{q-\frac 14} + \tfrac 1{\sqrt{c}} (t+T_c)^{\frac 12} 
 \sqrt{\tilde g_2(t)}\right) \\ &\leq
K \left(1+  \alpha c^{\frac 34 q-\frac 14} + \frac 1 {\sqrt c} \left(\int_0^{+\infty} \tilde g_2(t) dt\right)^{\frac 12}
+ \frac 1 {c} \left(\int_0^{+\infty} \tilde g_1(t) dt\right)^{\frac 12} + \exp(- c^{-r})\right).
\end{split}\end{equation*}
Thus,  from Lemma \ref{OdeT} and $-q+\frac 34 \geq 0$ ($p=2,3,4$) we obtain
$$
c^{-q+\frac 34} h_2(t)\leq K\left(1 + \alpha c^{-\frac q4}
+ c^{-q-\frac 14} \gamma_0^{\frac 12}(\alpha,c)+ \alpha^{-\frac 18}\exp(-\tfrac 12 c^{-r})\right).
$$
This estimate, inserted in \eqref{firstjj}, gives (i) for $J_2'$.


\begin{thebibliography}{10}
\bibitem{BL} H. Berestycki and P.-L. Lions,
Nonlinear scalar field equations. I. Existence of a ground state, Arch. Rational Mech. Anal. \textbf{ 82}, (1983) 313--345.
\bibitem{Craig} W. Craig, P. Guyenne, J. Hammack, D. Henderson and C. Sulem, Solitary water wave interactions.  Phys. Fluids  \textbf{18},  (2006),  057106.
\bibitem{Hirota}
R. Hirota, 
Exact solution of the Korteweg-de Vries equation for multiple collisions of solitons, Phys. Rev. Lett., \textbf{27} (1971), 1192--1194.
\bibitem{KATO}  T. Kato, On the Cauchy problem for the (generalized) Korteweg-de Vries equation. Studies in applied mathematics, 93--128, Adv. Math. Suppl. Stud., \textbf{8}, Academic Press, New York, 1983.
\bibitem{KPV} C.E. Kenig, G. Ponce and L. Vega, Well-posedness and scattering results for the generalized Korteweg--de Vries equation via the contraction principle. Comm. Pure Appl. Math. \textbf{46} (1993), 527--620.
\bibitem{Lax} P. D. Lax, Integrals of nonlinear equations of evolution and solitary waves, Comm. Pure Appl. Math. \textbf{21}, (1968) 467--490.
\bibitem{LiSattinger} Yi Li and D.H. Sattinger,  
Soliton collisions in the ion acoustic plasma equations. 
J. Math. Fluid Mech. \textbf{1}  (1999), 117--130. 
\bibitem{Martel} Y. Martel, Asymptotic $N$--soliton--like solutions of the subcritical and critical generalized Korteweg--de Vries equations, Amer. J. Math.  \textbf{127} (2005), 1103-1140.
\bibitem{yvanSIAM} Y. Martel, Linear problems related to asymptotic stability of solitons
of the generalized KdV equations, SIAM J. Math. Anal. \textbf{38} (2006), 759--781.
\bibitem{MMjmpa} Y. Martel and F. Merle, A Liouville theorem for the critical generalized Korteweg--de Vries equation. J. Math. Pures Appl. \textbf{79} (2000), 339--425.
\bibitem{MMarchives} Y. Martel and F. Merle, Asymptotic stability of solitons for subcritical generalized KdV equations. Arch. Ration. Mech. Anal. \textbf{157} (2001), 219--254.
\bibitem{MMgafa}  Y. Martel and F. Merle, Instability of solitons for the critical generalized Korteweg-de Vries equation. Geom. Funct. Anal. 11 (2001),  74--123.
\bibitem{MMjams} Y. Martel and F. Merle, Blow up in finite time and dynamics of blow up solutions for the  $L^2$--critical generalized KdV equations. J. Amer. Math. Soc. \textbf{15} (2002), 617--664.
\bibitem{MMnonlinearity} Y. Martel and F. Merle, Asymptotic stability of solitons of the subcritical gKdV equations revisited. Nonlinearity \textbf{18} (2005), no. 1, 55--80.
\bibitem{MMas1} Y. Martel and F. Merle, Asymptotic stability of solitons of the gKdV equations with a general nonlinearity. To appear in Math. Annalen.
\bibitem{MMcol1} Y. Martel and F. Merle, Description of two soliton collision for the
quartic gKdV equation, arXiv:0709.2672v1
\bibitem{MMcol2} Y. Martel and F. Merle, Stability of two soliton collision for   nonintegrable gKdV equations,  arXiv:0709.2677v1.
\bibitem{MMT} Y. Martel, F. Merle and Tai-Peng Tsai,
Stability and asymptotic stability in the energy space of the sum of $N$ solitons
for the subcritical gKdV equations, Commun. Math. Phys. \textbf{231}, (2002) 347--373.
\bibitem{Miura} R.M. Miura, The Korteweg--de Vries equation: a survey of results, SIAM Review \textbf{18}, (1976) 412--459.
\bibitem{Mizu} T. Mizumachi, Large time asymptotics of solutions around solitary waves to the generalized Korteweg-de Vries equations, SIAM J. Math. Anal. \textbf{32} (2001),  1050--1080. 
\bibitem{PW} R.L. Pego and M.I. Weinstein, Asymptotic stability of solitary waves. Commun. Math. Phys. \textbf{164} (1994), 305--349.
\bibitem{We1} M.I. Weinstein, Lyapunov stability of ground states of nonlinear dispersive evolution equations. Comm. Pure Appl. Math. \textbf{39} (1986), 51--68.
\end{thebibliography}
\end{document}